\documentclass{article}

\usepackage[top=1in,bottom=1in,left=1in,right=1in]{geometry}

\usepackage{graphicx}
\usepackage{amsthm}
\usepackage{amsfonts}
\usepackage{amssymb}
\usepackage{amsmath}
\usepackage{float}

\usepackage[font = small , width = 10.5cm]{caption}

\newcommand{\mycolon}{\!:\!}

\theoremstyle{definition}
\newtheorem{theorem}{Theorem}[section]
\newtheorem{lemma}[theorem]{Lemma}

\newtheorem{construction}[theorem]{Construction}

\newtheorem{remark}[theorem]{Remark}

\begin{document}

\title{Applications Of Ordinary Voltage Graph Theory To Graph Embeddability, Part 1}

\author{Steven Schluchter\\ Department of Mathematical Sciences\\ George Mason University\\ Fairfax, VA 22030\\ sschluch@gmu.edu.}

\maketitle

\abstract{We study embeddings of a graph $G$ in a surface $S$ by considering representatives of different classes of $H_1(S)$ and their intersections.  We construct a matrix invariant that can be used to detect homological invariance of elements of the cycle space of a cellularly embedded graph.  We show that: for each positive integer $n$, there is a graph embeddable in the torus such that there is a free $\mathbb{Z}_{2p}$-action on the graph that extends to a cellular automorphism of the torus; for an odd prime $p$ greater than 5 the Generalized Petersen Graphs of the form $GP(2p,2)$ do cellularly embed in the torus, but not in such a way that a free-action of a group on $GP(2p,2)$ extends to a cellular automorphism of the torus; the Generalized Petersen Graph $GP(6,2)$ does embed in the the torus such that a free-action of a group on $GP(6,2)$ extends to a cellular automorphism of the torus; and we show that for any odd $q$, the Generalized Petersen Graph $GP(2q,2)$ does embed in the Klein bottle in such a way that a free-action of a group on the graph extends to a cellular automorphism of the Klein bottle.\medskip

}

\section{Introduction}\label{section:introduction}


An ordinary voltage graph encodes a highly symmetric covering graph, called the derived graph of the ordinary voltage graph.  An ordinary voltage graph embedding in a surface $S$ encodes a another graph embedding, called the derived embedding, which produces a highly symmetric cellular decomposition of the derived surface, which is a type of highly-symmetric covering space of $S$.  One of the early results that followed from the formulation of the construction of a derived graph is Theorem \ref{theorem:VoltageGraphs}, which states that each free-action of a group on a graph can be encoded in the form of an ordinary graph.  Another one of the earliest formulated consequences of the construction is Theorem \ref{theorem:VoltageGraphsExtension}, which implies that if a cellular embedding of a graph in a surface has the property that a free-action of a group on the graph extends to a cellular automorphism of the surface, then the embedding can be encoded in the form of an ordinary voltage graph embedding.  In 1992, Archdeacon in \cite{A} used ordinary voltage graph constructions to find embeddings of complete bipartite graphs with predictable duals.  Recently, in \cite{AS1}, Abrams and Slilaty used voltage graph constructions to classify the cellular automorphisms of the surfaces of Euler characteristic at least $-1$, and in \cite{AS4}, they used other voltage graph constructions to find the minimal graphs with a $\mathbb{Z}_n$-symmetry that cannot be embedded in the sphere in such a way that the symmetry extends to a $\mathbb{Z}_n$-symmetry of the sphere.  In practice, the derived embedding of an ordinary voltage graph embedding is hard to understand through the readily available information: a rotation scheme or a list of facial boundary walks combinatorially encodes a unique embedding of the derived graph in the derived surface in a surface, but neither the global combinatorial structure of the derived graph nor the homeomorphism class of the derived surface is immediately transparent.  Cellular homology does provide a means by which one can determine global information about a cellular graph embedding $G\rightarrow S'$ in an obscured surface $S'$.  However, the computation $H_1(S')$ requires the local data about the incidence of all faces and edges of $G\rightarrow S'$.  In this article, we utilize a known homological invariant called the $\mathbb{Z}_2$-intersection product (developed in \cite{Dubrovin}) defined for a surface $S'$ that may be obscured.  This test, though it is based on a somewhat crude invariant, does provide for the use of a small amount of intersection data to show that a set of classes of $H_1(S')$ is independent, which in turn can be used to derive global information about $S'$.  In short, the test is an application of our Theorem \ref{theorem:MatrixIndependence}, which states that a set of homology classes is independent if (but not only if) a certain matrix has rows that are linearly independent over $\mathbb{Z}_2$.  We use this invariant extensively in order to organize large families of graph embeddings in a non-obscured surface according to which 1-chains represent which homology classes of that surface.  We also use this invariant to derive global information about an obscured surface through understanding local intersection properties; since the torus has first Betti number 2, we can determine that a surface $S$ is not the torus by showing that the first Betti number of $S$ is greater than 2.

In Section \ref{section:GraphsAndHomology}, we develop all necessary graph theory, topological graph theory, and homology theory, which includes our Theorem \ref{theorem:MatrixIndependence}.  In Section \ref{section:VoltageGraphs}, we develop the basics of ordinary voltage graph theory.  In Section \ref{section:applications}, we state and prove our Theorem \ref{theorem:ThePoint}, which contains the main results of this paper: for each positive integer $n$, there is a graph embeddable in the torus such that there is a free $\mathbb{Z}_n$-action on the graph that extends to a cellular automorphism of the torus; for an odd prime $p$ greater than 5 the Generalized Petersen Graphs of the form $GP(2p,2)$ do cellularly embed in the torus, but not in such a way that a free-action of a group on $GP(2p,2)$ extends to a cellular automorphism of the torus; the Generalized Petersen Graph $GP(6,2)$ does embed in the the torus such that a free-action of a group on $GP(6,2)$ extends to a cellular automorphism of the torus; and we show that for any odd $q$, the Generalized Petersen Graph $GP(2q,2)$ does embed in the Klein bottle in such a way that a free-action of a group on the graph extends to a cellular automorphism of the Klein bottle, the nonorientable surface with the same Euler Characteristic as the torus.

\section{Basic graph, graph embedding, and homology theory}\label{section:GraphsAndHomology}

\subsection{Graphs and graph embeddings} 

For the purposes of this article, a graph $G=(V,E)$ is a finite and connected multigraph.  An edge is a link if it is not a loop.  A path in $G$ is a subgraph of $G$ that can be described as a sequence of vertices and edges $v_1e_1v_2e_2\ldots e_{l-1}e_l$ such that the $v_i$ are all distinct and the edge $e_i$ connects $v_i$ and $v_{i+1}$.  Given a subgraph $H$ of $G$, an $H$-path is a path that meets $H$ at its end vertices and only at its end vertices.  Let $D(G)$ denote the set of darts (directed edges) on the edges of $G$.  To each dart is associated a head vertex $h(d)$ and a tail vertex $t(d)$; we say that two distinct darts on the same edge are opposites of eachother.  We will call one dart on $e$ the positive dart and the other the negative dart.  A walk $W$ in $G$ is a sequence of darts $d_1d_2\ldots d_m$ such that $h(d_b)=t(d_{b+1})$ for all $b\in [m-1]$.  If $h(d_m)=t(d_1)$, then we say that $W$ is a closed walk.  The notion of an $H$-walk is defined by analogy with the definition of an $H$-path.

Let $S$ denote a compact and connected surface without boundary.  For the purposes of this article, $S^2$ $P^2$, $T$ and $\mbox{\textit{KB}}$ shall denote the sphere, the projective plane, the torus, and the Klein bottle, respectively.  Given that a graph is a topological space (a 1-complex), then we may define a graph embedding in a surface to be a continuous injection $i\colon G\rightarrow S$.  A cellular embedding of $G$ in $S$ is an embedding that subdivides $S$ into 2-cells.  The regions of the complement of $i(G)$ in $S$ are called the faces of $i$.  We will let $G\rightarrow S$ denote a cellular embedding of $G$ in $S$, and we will let $F(G\rightarrow S)$ denote the set of faces of $G\rightarrow S$.  Since $G\rightarrow S$ induces a cellular decomposition of $S$, $S$ can be thought of as a cellular chain complex $\left( \{C_k(S),\partial\right \rbrace$, where $C_0(S)$, $C_1(S)$, and $C_2(S)$ are the $\mathbb{Z}_2$-vector spaces of formal linear combinations of elements of $V(G)$, $E(G)$, and $F(G\rightarrow S)$, respectively, and $\partial$ is the usual boundary map $\partial\colon C_k(S)\rightarrow C_{k-1}(S)$.  Note that the fact that we are using $\mathbb{Z}_2$ coefficients means that for $f\in C_2(S)$, $\partial f$ can be expressed as a sum of all edges appearing exactly once in a boundary walk of $f$, regardless of any orientation on the edges bounding $f$.  For $X\subset E(G)$ or $X\in C_1(G)$, let $G\mycolon X$ denote the induced subgraph of $G$ consisting of edges of $X$ and the vertices to which the edges are incident.  For $X\subset F(G\rightarrow S)$ or $X\in C_2(S)$, we let $S\mycolon X$ denote the sub 2-complex consisting of the faces appearing in $X$, and all of their subfaces.

Given $G$, an automorphism of $G$ is a map $\phi\colon G\rightarrow G$ that bijectively maps $C_0(G)$ and $C_1(G)$ to $C_0(G)$ and $C_1(G)$, respectively, such that the incidence of edges and vertices is preserved.  Given $G\rightarrow S$ and considering the resulting 2-complex, a cellular automorphism of $S$ is a homeomorphism of $S$ that is an automorphism of $G$ and a bijection on $C_2(S)$ that preserves the incidence of faces of $G\rightarrow S$ with edges and vertices of $G$. A group $A$ is said to act cellularly on $S$ if there exists a graph $G'$ embedded in $S$ such that an action of $A$ on $G$ extends to a cellular automorphism of $S$.    

A subgraph of $G$ is called a \textit{circle} if it is a connected 2-regular graph.  Let $Z(G)$ denote the subsapce of $C_1(G)$ with generating set $\left \{z_i: G\mycolon z_i\ \mbox{is a circle in }G\right \rbrace$.  The subspace $Z(G)$ is called the cycle space of $G$.  For a circle $G\mycolon z$ of a cellularly embedded $G$, we let a ribbon neighborhood of $z$, denoted $R(z)$ denote a regular neighborhood of $G\mycolon z$ containing $G\mycolon z$, the ends of edges incident to the vertices of $G\mycolon z$, and no other edge segments or vertices.  Clearly, a ribbon neighborhood of any circle is homeomorphic to an annulus or a M\"obius band.  If a circle $C$ has a ribbon neighborhood that is homeomorphic to a M\"obius band, then $C$ is an orientation-reversing circle.  Else, $C$ is called an orientation-preserving circle.

Given $G\rightarrow S$ and a dart $d$ with $t(d)=v$, let $\rho\colon D\rightarrow D$ be the permutation that takes $d$ to the next dart in a cyclic order of darts with tail vertex $v$.  The order of darts with tail vertex $v$ that follows the order induced by $\rho$ is called the rotation on $v$.  A rotation on a vertex $v$ is a list of the form $v: d_id_jd_k\ldots $.  The permutation $\rho$ combinatorially encodes $G\rightarrow S$ and is called a rotation scheme on $G$.  Following \cite[\S 3.2]{GT}, if one thickens the embedded graph such that the vertices become discs and the edges become rectangular strips glued to the discs, one produces what is called a band decomposition of the surface $S$.  The 0-bands are the discs, the 1-bands are the rectangular strips, and the 2-bands are discs glued to the 1-bands and the 0-bands.  If one of the two possible orientations on any given 1-band is consistent with the orientations induced by $\rho$ on the joined 0-bands (see \cite[Figures 3.13 and 3.14]{GT} for enlightening diagrams), then it is said that the given 1-band is orientation preserving, else it is orientation reversing.  Thus, an edge $e$ of an embedded graph may be designated as an orientation-reversing edge if its associated 1-band is orientation reversing, else it is an orientation-preserving edge.  It may help to think of an orientation-reversing edge as being ``twisted".  If $e$ is an orientation-reversing edge of $G\rightarrow S$, then the rotation scheme of $G$ features a $1$ superscript above all occurrences of a dart on $e$.

\subsection{Homology theory and intersection theory}\label{subsection:HomologyAndIntersectionTheory}

Consider $G\rightarrow S$ and the associated 2-complex $\left \{C_k(S),\partial\right \rbrace$.  Let $B(G\rightarrow S)$ be the subspace of $Z(G)$ with generating set $\left \{\partial f:\ f\in C_2(S)\right \rbrace$.  We let $H_1(S)$ denote the first homology group of $S$.  Using our notation, $H_1(S) = Z(G) / B(G\rightarrow S)$.  Let $\beta_1(G)$ denote the rank of $H_1(S)$; $\beta_1(S)$ is commonly called the first Betti number of $S$.  We let $[\cdot]\colon Z(G)\rightarrow Z(G) / B(G\rightarrow S)$ denote the map that takes each element of $Z(G)$ to its equivalence class modulo $B(G\rightarrow S)$.  If $z\in Z(G)$ satisfies $[z]=[0]$, then we say that $z$ is homologically trivial.  If $z\in Z(G)$ is such that $G\mycolon z$ is a circle and $S\setminus (G\mycolon z)$ is connected, we say that $G\mycolon z$ is a nonseparating circle and a separating circle $S\setminus (G\mycolon z)$ is connected.  Recall that if $z\in Z(G)$ has the property that $G\mycolon z$ is a separating circle than $z$ is a homologically trivial and homologically nontrivial if $G\mycolon z$ is a nonseparating circle.  If $z_1,\ z_2\in Z(G)$ satisfy $[z_1]\neq[z_2]$, then we say that $z_1$ and $z_2$ are homologically independent, else we say that $z_1$ and $z_2$ are homologous.

As described in \cite{Dubrovin2}, for $z_1, z_2\in Z(G)$ inducing circles in $Z(G)$ which are embedded in general position in $S$, the intersection index $(G\mycolon z_1,G\mycolon z_2)$ is invariant under homotopy, and is the residue modulo 2 of the number of points (vertices) in $G\mycolon z_1\cap G\mycolon z_2$.

\begin{figure}[H]
\begin{center}
\includegraphics[scale=.5]{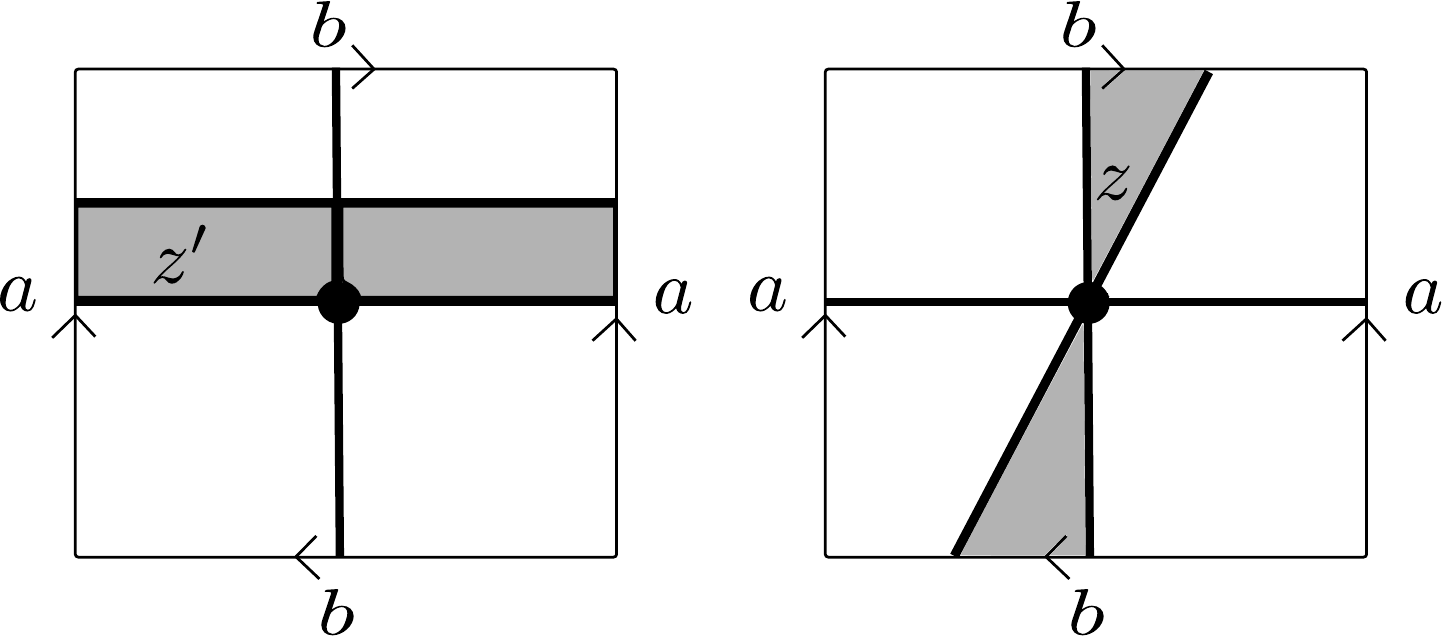}
\caption{A depiction illustrating the conclusions of Remark \ref{remark:TransverseCrossings}}.\label{fig:MobiusNeighborhood}
\end{center}
\end{figure}

\begin{remark}\label{remark:TransverseCrossings} Regardless of orientability of $S$, as evidenced in Figure \ref{fig:MobiusNeighborhood}, a circle $G\mycolon z$ of $G\rightarrow S$ is orientation reversing iff $(G\mycolon z,G\mycolon z)=1$, and a circle $G\mycolon z'$ orientation preserving iff $(G\mycolon z',G\mycolon z')=0$.\end{remark}

Also discussed in \cite[p. 85]{Dubrovin2}, any two simple closed curves that are embedded in a surface can be brought into general position with each other by using a homotopy of one of the embeddings, and it is a consequence of the theory that this does not affect the intersection index of the two curves.  The intersection product can therefore be defined for any two simple closed curves in a surface.  Since $H_1(S)$ is isomorphic to the abelianization of the fundamental group $\pi_1(S)$, we see that the intersection index of two circles of $G$ can also be extended to 1-chains inducing circles, treating them as representatives of classes of $H_1(S)$.  Per \cite[Theorem 18.1]{Dubrovin} and the discussions that follow on \cite[pp.220-221]{Dubrovin}, after making a choice of basis of $H_1(S)$, this intersection index is a well defined bilinear map giving a nondegenerate, symmetric, bilinear form on $H_1(S)$, $\langle \cdot, \cdot \rangle\colon H_1(S)\times H_1(S)\rightarrow \mathbb{Z}_2$.  For the sake of brevity, we will let $\langle z_1,z_2\rangle$ denote $\langle [z_1],[z_2]\rangle$ for the remainder of this article.

\begin{lemma}\label{lemma:ZIntersections}  Consider $G\rightarrow S$. If $z'\in Z(G)$ is homologically trivial, then $\langle z' ,z\rangle=0$ for all $z\in Z(G)$.\end{lemma}

\proof This follows immediately from the bilinearity of the intersection indices since $[z']=[0]$.\endproof\medskip

For $G\rightarrow S$ let $\left \{x_1,\ldots,x_l\right \rbrace \subset Z(G)$ denote an ordered set of 1-chains, with each element $x_i\in X$ inducing a circle in $G$.  Let $M_X$ denote the symmetric $l\times l$ matrix over $\mathbb{Z}$ whose entry $m_{ij}$ is $\langle x_i,x_j\rangle$.\medskip

\begin{theorem}\label{theorem:MatrixIndependence}If the rows of $M_X$ are linearly independent, then the elements of $X$ are homologically independent.\end{theorem}

\proof  Assume for the sake of contradiction that the nontrivial $\mathbb{Z}_2$-sum $\sigma = \sum a_i x_i$ is homologically trivial. Write $\sigma =(a_1,a_2,\ldots,a_l) \in \mathbb{Z}_2^l$.  Since the rows of $M_X$ are linearly independent, then $M_X \sigma^T$ is a nonzero vector in $\mathbb{Z}_2^l$.  Thus, there exists a nonzero $y=(y_1,y_2,\ldots, y_l) \in \mathbb{Z}_2^l$ such that \[y M_X \sigma^T\neq 0.\]  Let $\sigma_y =\sum y_ix_i$.  It follows that \begin{align*} y  (M_X \sigma) & =  \sum_{i,j}a_iy_j (G\mycolon x_i , G\mycolon x_j)\\ &=  \langle \sigma_y,\sigma \rangle \\ & \neq 0.\end{align*}  This contradicts Lemma \ref{lemma:ZIntersections}.\endproof\medskip \endproof

\begin{lemma}\label{lemma:P2CirclesCross}\cite[p. 222, paragraph (c)]{Dubrovin} Consider $G\rightarrow P^2$ and $z_1, z_2\in Z(G)$.  If $z_1$ and $z_2$ both induce orientation-reversing circles, then $\langle z_1,z_2\rangle=1$.\end{lemma}

\begin{lemma}\label{lemma:TorusCirclesCross}Consider $G\rightarrow T$.  If $z_1$ and $z_2$ are homologically nontrivial, homologically independent, and induce circles, then $\langle z_1,z_2\rangle =1$.\end{lemma}

\proof Assume for the sake of contradiction that $\langle z_1,z_2\rangle=0$.  Since $z_1$ is homologically nontrivial and orientation-preserving, $G\mycolon z_1$ is a nonseparating curve and so there exists $z_3\in Z(G)$ such that $G\mycolon z_3$ is a circle and $\langle z_3,z_1\rangle=1$, yet we can assume nothing about the value of $\langle z_2,z_3\rangle$.  Similarly, there is a $z_4\in Z(G)$ such that $G\mycolon z_4$ is a circle and $\langle z_2,z_4\rangle = 1$, yet we can assume nothing about $\langle z_1,z_4\rangle$ and $\langle z_3,z_4\rangle$.  Let $X=\left \{z_1,z_2,z_3,z_4\right \rbrace$ and note that $M_{X}$ is a symmetric matrix over $\mathbb{Z}_2$ of the form \[ \begin{pmatrix} 0 & 0 & 1 & * \\ 0 & 0 & * & 1 \\ 1 & * & 0 & * \\ * & 1 & * & 0 \end{pmatrix},\] where the elements of $M_{X}$ denoted by $*$ are unknown elements of $\mathbb{Z}_2$.

Since the rows of this matrix are linearly independent, Theorem \ref{theorem:MatrixIndependence} implies that $\beta_1(T) \ge 4$, which contradicts the fact that $\beta_1(T)=2$.\endproof

\begin{lemma}\label{lemma:KBCirclesCross} Consider $G\rightarrow \mbox{KB}$ and $z_1,z_2\in Z(G)$.  If $z_1$ induces an orientation-reversing circle and $z_2$ is homologically nontrivial and induces an orientation-preserving circle, then $\langle z_1,z_2\rangle=1$.\end{lemma}

\proof  Assume for the sake of contradiction that $\langle z_1,z_2\rangle=0$.  By Remark \ref{remark:TransverseCrossings}, $\langle z_1,z_1\rangle=1$ and $\langle z_2,z_2\rangle=0$.  Since $z_2$ is homologically nontrivial, $G\mycolon z_2$ is a nonseparating curve and so there exists $z_3\in Z(G)$ such that $G\mycolon z_3$ is a circle and $\langle z_2,z_3\rangle=1$, yet we can assume nothing about the values of $\langle z_1,z_3\rangle$ and $\langle z_3,z_3\rangle$.  Let $X=\left \{z_1,z_2,z_3\right \rbrace$ and note that $M_{X}$ is a symmetric matrix over $\mathbb{Z}_2$ of the form \[ \begin{pmatrix} 1 & 0 & *\\ 0 & 0 & 1 \\ * & 1 & * \end{pmatrix},\] where the elements of $M_{X}$ denoted by $*$ are unknown elements of $\mathbb{Z}_2$.

Since the rows of this matrix are linearly independent, Theorem \ref{theorem:MatrixIndependence} implies that $\mbox{dim}_{\mathbb{Z}_2}(H_1\mbox{\it KB})\ge 3$, which contradicts the fact that $\mbox{dim}_{\mathbb{Z}_2}(H_1(\mbox{\it KB}))=2$.\endproof

For a cellular complex $K$, let $\chi(K)$ denote the Euler characteristic of $K$, which is defined \[\chi(K) = \sum_{i}(-1)^i k_i,\] where $k_i$ is the number of cells of dimension $i$ in $K$.  We let $\beta_i(K)$ denote the $i^{th}$ Betti number of $K$. Theorem \ref{theorem:EulerCharConnectedSum} can be deduced from the definition of Euler characteristic and the well-known result \cite[Theorem 2.4.4]{H}, which says that \[\chi(K)=\sum_{i} (-1)^i\beta_i(K).\]

\begin{theorem}\label{theorem:EulerCharConnectedSum}For a connected surface $S$ with $h$ handles and $c$ crosscaps, \[\chi(S)=2-2h-c.\]\end{theorem}

Another tool that we will make use of is the well-known Classification of Surfaces, with and without boundary, which can be found in \cite[Theorem 5.1]{M}.  Lemmas \ref{lemma:P2Decomposition}, \ref{lemma:KBDecomposition}, and \ref{lemma:TDecomposition}, which concern the decomposition of a surface $S$ as a connected sum $S_1\#S_2$ of two surfaces $S_1$ and $S_2$, are consequences of the Classification of Surfaces, Theorem \ref{theorem:EulerCharConnectedSum}, and the (non)orientability of the projective plane, the torus, and the Klein bottle.

\begin{lemma}\label{lemma:P2Decomposition}As a connected sum of surfaces, the projective plane $P^2$ is only decomposable as $P^2\#S^2$.\end{lemma}

\begin{lemma}\label{lemma:KBDecomposition}As a connected sum of surfaces, the Klein bottle $\mbox{\textit{KB}}$ has exactly two decompositions: $\mbox{\textit{KB}}=\mbox{\textit{KB}}\#S^2$ and $\mbox{\textit{KB}}=P^2\#P^2$.\end{lemma}

\begin{lemma}\label{lemma:TDecomposition}As a connected sum of surfaces, the torus $T$ is only decomposable as $T=T\#S^2$.\end{lemma}

\section{Basic ordinary voltage graph theory}\label{section:VoltageGraphs}

Consider $G$ and let $e$ denote an edge of $G$.  Following $\cite{GT}$, we let $e$ denote the positive edge on $e$ and $e^{-}$ denote the negative edge on $e$.  Let $A$ denote a finite group and let $1_A$ denote the identity element of $A$.  An ordinary voltage graph is an ordered pair $\langle G, \alpha \rightarrow A \rangle$ such that $\alpha\colon D\rightarrow A$ satisfies $\alpha(e^-)=\alpha(e)^{-1}$.  The group element $\alpha(e)$ is called the voltage of $e$.  Associated to each ordinary voltage graph is a derived graph $G^\alpha=(V\times A, E\times A)$.  The directed edge $(e,a)$ has tail vertex $(v,a)$ and head vertex $(v,a\alpha(e))$; as a consequence of this and the conditions imposed on $\alpha$, the dart $(e^{-},a\alpha(e)^{-1})$ is the dart opposite $(e,a)$.  We will use the abbreviation $v^a$ for $(v,a)$ and $e^a$ for $(e,a)$.  We let $p\colon G^\alpha \rightarrow S$ denote the projection (covering) map satisfying $p(e^a)=e$ and $p(v^a)=v$.  For a walk $W=d_1d_2\ldots d_m$, let $\omega(W)=\alpha(d_1)\alpha(d_2)\ldots\alpha(d_m)$ denote the net voltage of $W$.  If $c\in C_1(G)$ is such that $G\mycolon c$ is connected and $W$ is a closed walk of $G\mycolon c$ based at $v\in V(G\mycolon c)$, then, per \cite[Theorem 2.1.1]{GT}, each lift of $W$ is uniquely identifiable by the vertex $v^a$ at which it begins.  For each $a\in A$, let $W_v^a$ denote the lift of $W$ beginning at $v^a$.  Note that $W_v^a$ ends at the vertex at which $W_v^{a\omega(W)}$ begins.

Following \cite{ASc}, we say that lifts $W_v^a$ and $W_v^{b}$ are consecutive if $a\omega(W)=b$ or $b\omega(W)=a$.  We call a set of lifts of $W$ of the form \[\left \{W_v^a,W_v^{a\omega(W)},W_v^{a(\omega(W))^2},\ldots\right \rbrace\] a set of consecutive lifts of $W$. Observe that $W_v^a$ and $W_v^b$ are lifts of $W$ contained in the same set of consecutive lifts of $W$ if and only if $a=b(\omega(W))^m$ for some nonnegative integer $m$.  Let $\hat{W}_v^{a}$ denote the set of consecutive lifts of $W$ containing $W_v^{a}$.  Since $A$ is assumed to be finite, each set of consecutive lifts of $W$ contains a finite number of lifts of $W$.  Unless $\omega(W)=1_A$, there is more than one group-element superscript that can be used to identify $\hat{W}_v^{a}$.  This conclusion also holds for all other sets of consecutive lifts of $W$.

Also described in \cite{GT}, an ordinary voltage graph embedding of $G$ in $S$ is an ordered pair $\langle G \rightarrow S, \alpha \rightarrow A\rangle$, which is referred to as a base embedding.  Each base embedding encodes a derived embedding, denoted $G^\alpha \rightarrow S^\alpha$, in the derived surface $S^\alpha$.  Gross and Tucker in \cite{GT} describe the derived embedding according to rotation schemes, but we use Garman's manner of describing it.  Garman points out in \cite{G} that since it is the lifts of facial boundaries that form facial boundaries in $S^\alpha$, $S^\alpha$ can be formed by ``idendifying each component of a lifted region with sides of a 2-cell (unique to that component) and then performing the standard identification of edges from surface topology".  It is therefore permissible to have a base embedding in a surface $\hat{S}$ with or without boundary; for each (directed) edge $e$ bounded on only one side by a face of $G\rightarrow \hat{S}$, each (directed) edge $e^a$ is bounded on only one side by a face of $G^\alpha \rightarrow \hat{S}^\alpha$.

Lemma \ref{lemma:OrientationReversingVoltage} is a special case of \cite[Theorem 4.1.4]{GT}.  \begin{lemma}\label{lemma:OrientationReversingVoltage}Consider an ordinary voltage graph embedding in a nonorientable surface.  If there exists a $z\in Z(G)$ such that $G\mycolon z$ is an orientation-reversing circle traversable by an Eulerian walk $W$ such that $|\langle \omega(W)\rangle |$ is odd, then the derived surface is nonorientable.\end{lemma}

For $\langle G\rightarrow S,\alpha \rightarrow A\rangle$ we let $S_v^a$ denote the component of $S^\alpha$ containing the vertex $v^a$, and for $\langle G,\alpha \rightarrow A\rangle$ we let $G_v^a$ denote the component of $G^\alpha$ containing $v^a$. We use similar notation for induced ordinary voltage graphs and ordinary voltage graph embeddings: for $I\in C_2(S)$ and $v\in V(S\mycolon I)$, $(S\mycolon I)_v^a$ is the component of $(S\mycolon I)^\alpha$ containing $v^a$.  For shorthand, if $x\in C_1(G)$, and $v\in V(G\mycolon x)$, we will use  $x_v^a$ to denote the 1-chain inducing $(G\mycolon x)_v^a$.

Following \cite{GT}, the voltage group $A$ acts by left multiplication on $G^\alpha$; for $c\in A$, let $c\cdot v^a=v^{ca}$, $c\cdot e^a= c\cdot e^{ca}$.  This group action is clearly regular (free and transitive) on the fibers over vertices and (directed) edges of $G^\alpha$, and so the components of $G^\alpha$ are isomorphic.  This action extends to a transitive (not necessarily free) action on the faces forming the fiber over a face of a base embedding, and so the components of $S^\alpha$ are homeomorphic as topological spaces and isomorphic as cellular complexes.  It is a consequence of the theory that the graph covering map can be extended to a (branched) covering map of surfaces
\cite[Corollary to Theorem 4.3.2]{GT}.  We will use $p\colon S^\alpha \rightarrow S$ to denote the associated covering map.

Consider an ordinary voltage graph embedding $\langle G\rightarrow S, \alpha \rightarrow A\rangle$.  For a fixed $v\in V(G)$, $A(v)$ denotes the local voltage group of net voltages of closed walks in $G$ based at $v$; for $I\in C_2(S)$ such that $S\mycolon I$ is connected and contains $v$, we let $A(v,I)$ denote the local voltage group of closed walks in $S\mycolon I$ based at $v$; for $y\in C_1(S\mycolon I)$ such that $G\mycolon y$ is connected and contains $v$, we let $A(v,y)$ denote the local voltage group of net voltages of closed walks in $G\mycolon y$ based at $v$.

\begin{theorem}\label{theorem:BigCosetTheorem}\cite[Theorem 3.8]{ASc} Consider $\langle G\rightarrow S, \alpha \rightarrow A\rangle$, and let $I$ be a subset of $F(G\rightarrow S)$ such that $S\mycolon I$ is connected.  Let $x\in C_1(S\mycolon I)$ induce a connected subgraph of $G$, $v$ denote a vertex of $G\mycolon x$, and $W$ be a closed walk of $G\mycolon x$ based at $v$.

\begin{enumerate}

\item{There are $\frac{|A|}{|A(v)|}$ components of $S^\alpha$.}\label{theorem:BigCosetTheoremPartOne}

\item{There are $\frac{|A|}{|A(v,I)|}$ components of $(S\mycolon I)^\alpha$ contained in each component of $P^\alpha$.}\label{theorem:BigCosetTheoremPartTwo}

\item{There are $\frac{|A(v,I)|}{|A(v,x)|}$ components of $(G\mycolon x)^\alpha$ contained in each component of $(S\mycolon I)^\alpha$.}\label{theorem:BigCosetTheoremPartThree}

\item{There are $\frac{|A(v,x)|}{|\langle \omega(W)\rangle|}$ sets of consecutive lifts of $W$ covering the edges of each component of $(G\mycolon y)^\alpha$.}\label{theorem:BigCosetTheoremPart4}
\end{enumerate} \end{theorem}

Before we state and prove Theorems \ref{theorem:ThePoint}, we state the remaining necessary background information.  As alluded to before, two of the central accomplishments of ordinary voltage graph theory are Theorems \ref{theorem:VoltageGraphs} and \ref{theorem:VoltageGraphsExtension}.

\begin{theorem}\label{theorem:VoltageGraphs}\cite[Theorem 2.2.2]{GT} If $A$ is a group acting freely on a graph $\tilde{G}$, and $G$ is the resulting quotient graph, then there is an assignment $\alpha$ of ordinary voltages in $A$ to the quotient graph $G$ such that $G^\alpha$ is isomorphic to $\tilde{G}$ and that the given action of $A$ on $\tilde{G}$ is the natural left action of $A$ on $G^\alpha$.\end{theorem}

For any derived embedding of an ordinary voltage graph embedding $\langle G\rightarrow S, \alpha \rightarrow A\rangle$, there is at most one branch point contained in each face of the base embedding.  Let $f_y$ denote a face of $G\rightarrow S$ containing a branch point $y$, $W_{f_y}$ denote a facial boundary walk of $f_y$ and $\omega_{f_y}$ denote the net voltage of $W_{f_y}$.  Let $|\omega{f_y}|$ denote $|\langle \omega_{f_y}\rangle |$, and recall that per Theorem \ref{theorem:BigCosetTheorem} that it takes $|\omega_{f_y}|$ lifts of $W_{f_y}$ to form a facial boundary of a face in the fiber over $f_y$.  Note that in $S^\alpha$, which is an $|A|$-fold branched covering of $S$,\[|A|-\frac{|A|}{|\omega_{f_y}|}=|p^{-1}(f_y)|.\] We define the \textit{deficiency} of $y$, which we denote $\mbox{\textit{def}(y)}$ to be $|A|-p^{-1}(f_y)$.  So, for an $n$-fold cover of $S$, \[\mbox{\textit{def}(y)}=n-|p^{-1}(y)|.\]

\begin{theorem}\label{theorem:RHEquation} (The Ruemann-Hurwicz Equation) \cite[Theorem 4.2.3]{GT} Let $p\colon \tilde{S}\rightarrow \dot{S}$ denote an $n$-fold branched covering of surfaces and let $Y$ denote the set of branch points of $\dot{S}$. Then \[\chi(\tilde{S})=n\chi(\dot{S})-\sum_{y\in Y}\mbox{\textit{def}(y)}.\]\end{theorem}

As discussed in \cite[\S 4.3]{GT}, a finite-sheeted covering space $(\tilde{S},p)$ of a surface $S$ is said to be regular if $S$ is the quotient of $\tilde{S}$ modulo the action of a finite group $A$ has a finite number of points.  Such an action on a surface is called \textit{pseudofree}.  Theorem \ref{theorem:VoltageGraphsExtension} establishes that under certain circumstances, a pseudofree-action of group on a surface can be encoded in the form of an ordinary voltage graph embedding.

\begin{theorem}\cite[Theorem 4.3.5]{GT}\label{theorem:VoltageGraphsExtension} Let $p\colon \tilde{S}\rightarrow S$ denote a regular branched covering of surfaces, and let $G\rightarrow S$ denote an embedding in $S$ having at most one branch point in any face and no branch points in $G$.  There is an ordinary voltage graph assignment $\alpha$ in the group of (covering) deck transformations of $p$ such that $p$ is equivalent to a natural surface projection $S^\alpha\rightarrow S$.\end{theorem} 

Applying Theorems \ref{theorem:VoltageGraphsExtension}, we see that if $G$ is embedded in $S$ in such a way that the free-action of a group $A$ on $G$ extends to a cellular automorphism of $S$, then the embedding can be encoded using an ordinary voltage graph embedding.

\section{Applications to graph embeddability}\label{section:applications}

\begin{figure}[H]
\begin{center}
\includegraphics[scale=.5]{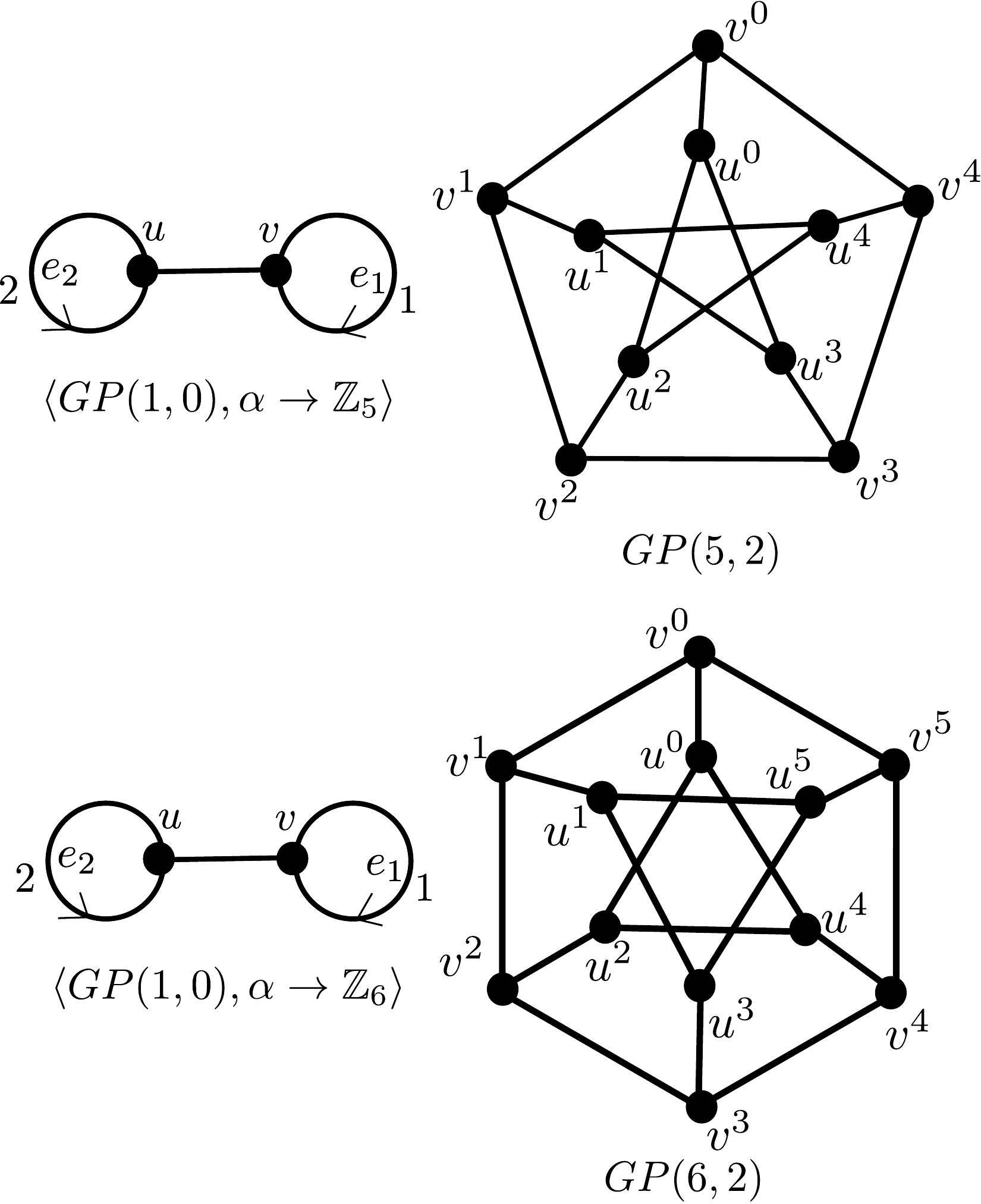}
\caption{The Petersen Graph and the D\"{u}rer Graph as derived graphs of ordinary voltage graphs.  The darts on the edges of each derived graph are omitted, as are the zero voltages in the base embedding.}\label{fig:GPExamples}
\end{center}
\end{figure}

Given positive integers $n,\ k$, the Generalized Petersen Graph $GP(n,k)$ has vertices $v^0,\ v^1,\ \ldots,\ v^{n-1}$, $u^0,\ u^1,\ \ldots, u^{n-1}$.  For each $i\in \mathbb{Z}_n$, $GP(n,k)$ has edges $(v^i,v^{i+1})$, $(v^i,u^i)$, and $(u^i,u^{i+k})$, where all superscripts are taken modulo $n$.  We use superscripts for these graphs to have a more unified set of notation so that we may work within our stated framework for ordinary voltage graphs.  Two well known Generalized Petersen Graphs are the Petersen Graph, $GP(5,2)$, and the D\"{u}rer Graph, $GP(6,2)$.  Both of these graphs appear as derived graphs in Figure \ref{fig:GPExamples}.

\begin{figure}[H]
\begin{center}
\includegraphics[scale=.5]{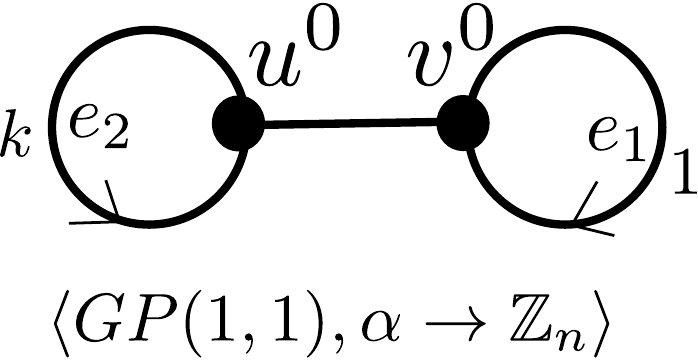}
\caption{An ordinary voltage graph whose derived graph is $GP(n,k)$.  The zero voltages are not shown.}\label{fig:GPBase}
\end{center}
\end{figure}

\begin{remark}\label{remark:ZnAction}Note the presence of a free $\mathbb{Z}_n$-action on $GP(n,k)$.  An integer $l$ acts on $GP(n,k)$ by mapping $v^i$ to $v^{i+l}$ and $u^i$ to $u^{i+l}$, where the superscripts are reduced modulo $n$.  Note that this action is transitive on the vertices $v^i$ and it is transitive on the vertices $u^i$.  The quotient of $GP(n,k)$ modulo this $\mathbb{Z}_n$-action is the ``barbell graph" $GP(1,0)$.  Each Generalized Petersen Graph can be recovered as the derived graph of the ordinary voltage graph appearing in Figure \ref{fig:GPBase}.\end{remark}

The results of this section are applications of ordinary voltage graph theory that are directed toward proving Theorem \ref{theorem:ThePoint}.

\begin{theorem}\label{theorem:ThePoint}Let $p$ and $q$ be odd primes.
     \begin{enumerate}
     \item{The group $\mathbb{Z}_{2p}$ acts cellularly on the torus.}\label{theorem:ThePointPart1}
     \item{Each Generalized Petersen Graph $GP(2p,2)$ has a cellular embedding in the torus.}\label{theorem:ThePointPart2}
     \item{For $p=3$, the Generalized Petersen Graph $GP(2p,2)$ has a cellular embedding in the torus in such a way that a free-action of a group on $GP(2p,2)$ extends to a cellular automorphism of the torus.}\label{theorem:ThePointPart3}
     \item{For $p>5$, the Generalized Petersen Graph $GP(2p,2)$ has no cellular embedding in the torus in such a way that a free-action of a group on $GP(2p,2)$ extends to a cellular automorphism of the torus.}\label{theorem:ThePointPart4}
      \item{Each Generalized Petersen Graph $GP(2p,2)$ has a cellular embedding in the Klein bottle in such a way that a free action of a group on $GP(2p,2)$ extends to a cellular automorphism of the Klein bottle.}\label{theorem:ThePointPart5}
      \end{enumerate}
\end{theorem}

\proof

\ \medskip \\ \textit{Proof of Part \ref{theorem:ThePointPart1}.}\medskip 

\begin{figure}[H]
\begin{center}
\includegraphics[scale=.45]{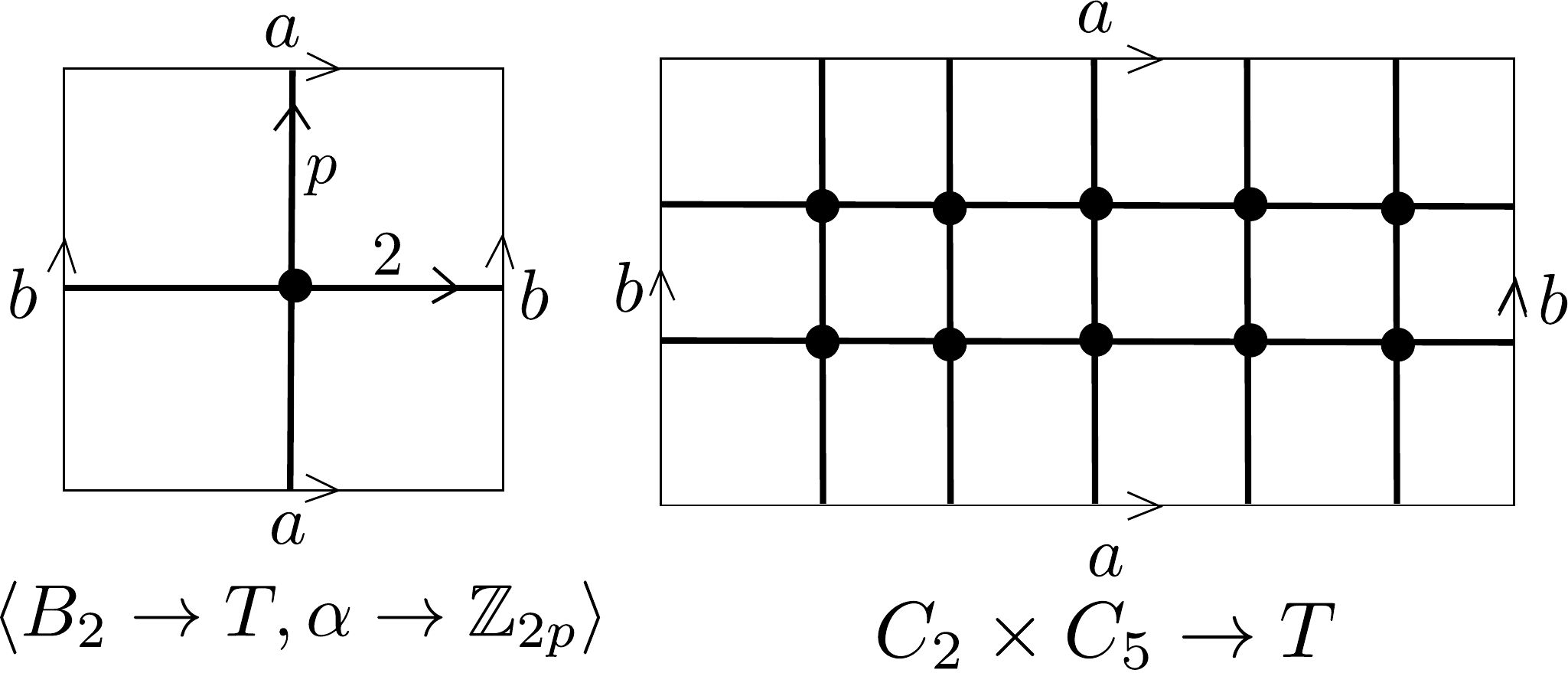}
\caption{An ordinary voltage graph embedding of the bouquet of two loops in the torus with voltage group $\mathbb{Z}_{2p}$ such that the derived embedding is in the torus.  The darts on the edges of the derived graph are omitted.  The derived embedding is the derived embedding of the special case of the base embedding for which $p=5$.}\label{fig:GPBase4}
\end{center}
\end{figure}

Let $p$ be an odd prime. Figure \ref{fig:GPBase4} proves Part \ref{theorem:ThePointPart1} of Theorem \ref{theorem:ThePoint}.

\noindent \textit{Proof of Part \ref{theorem:ThePointPart2}.}\medskip  

\begin{figure}[H]
\begin{center}
\includegraphics[scale=.38]{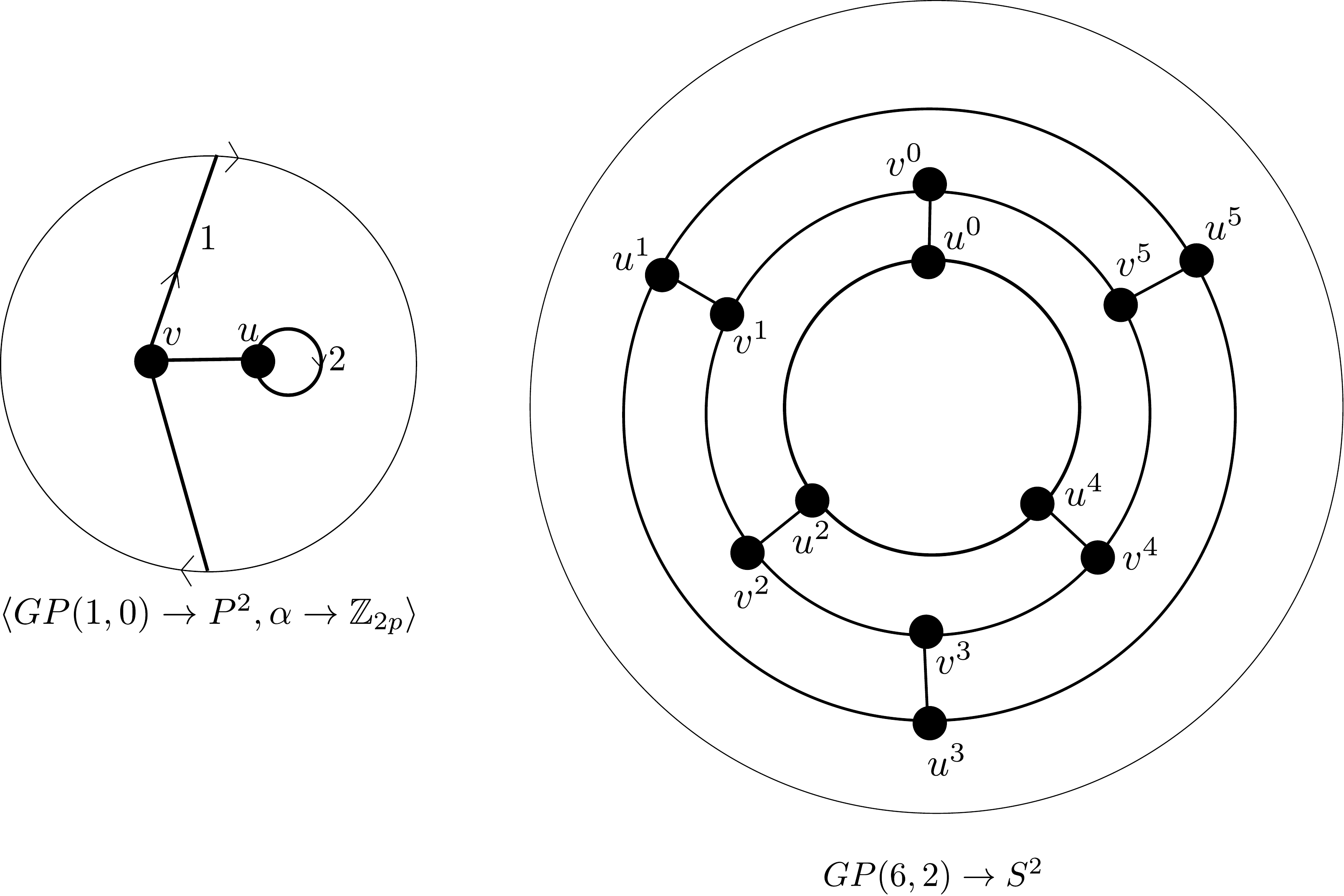}
\caption{$GP(2p,2)$ embedded in the sphere as a derived embedding of the special case of the base embedding for which $p=3$.  The zero voltages are not shown.  The darts on the edges of the derived graph are omitted.}\label{fig:SphereDerived}
\end{center}
\end{figure}

We first show that that $GP(2p,2)$ has a derived embedding in the sphere as evidenced by Figure \ref{fig:SphereDerived}.  In Construction \ref{constructon:GPTorus}, we will alter the derived embeddings appearing in Figure \ref{fig:SphereDerived} to show that for each odd prime $p$, $GP(2p,2)$ has a cellular embedding in the torus.

\begin{remark}\label{remark:SphereDerived}The two faces of the base embedding in Figure \ref{fig:SphereDerived} have boundary walks whose voltages generate the subgroups $\langle 0 \rangle $ and $\langle 2\rangle$ of $\mathbb{Z}_{2p}$.  By Theorem \ref{theorem:RHEquation}, we see that in this case, $(P^2)^\alpha$ has Euler characteristic $2p\cdot 1 - (2p-2)=2$, and so, $(P^2)^\alpha = S^2$ by Theorem \ref{theorem:EulerCharConnectedSum}.\end{remark}

Using an adaptation of a technique of Xuong \cite[Theorem 3.4.13]{GT}, we now show that each $GP(2p,2)$ has a cellular embedding in the torus.

\begin{construction}\label{constructon:GPTorus}  Given the derived embedding of $GP(2p,2)$ in the sphere described in Remark \ref{remark:SphereDerived}, we construct a cellular embedding of $GP(2p,2)$ in the torus.  Form the connected sum $S^2\# T$ with the requirement that the circle of attachment is contained in the face $f$ bounded by the circle of $GP(2p,2)$ that contains the vertices $u^{2l}$.  Now remove the edges $(u^0,u^2)$, $(u^2,u^4)$, as in Figure \ref{fig:GPStep2}; note that the face $f$ is homeomorphic to a punctured torus.  

\begin{figure}[H]
\begin{center}
\includegraphics[scale=.4]{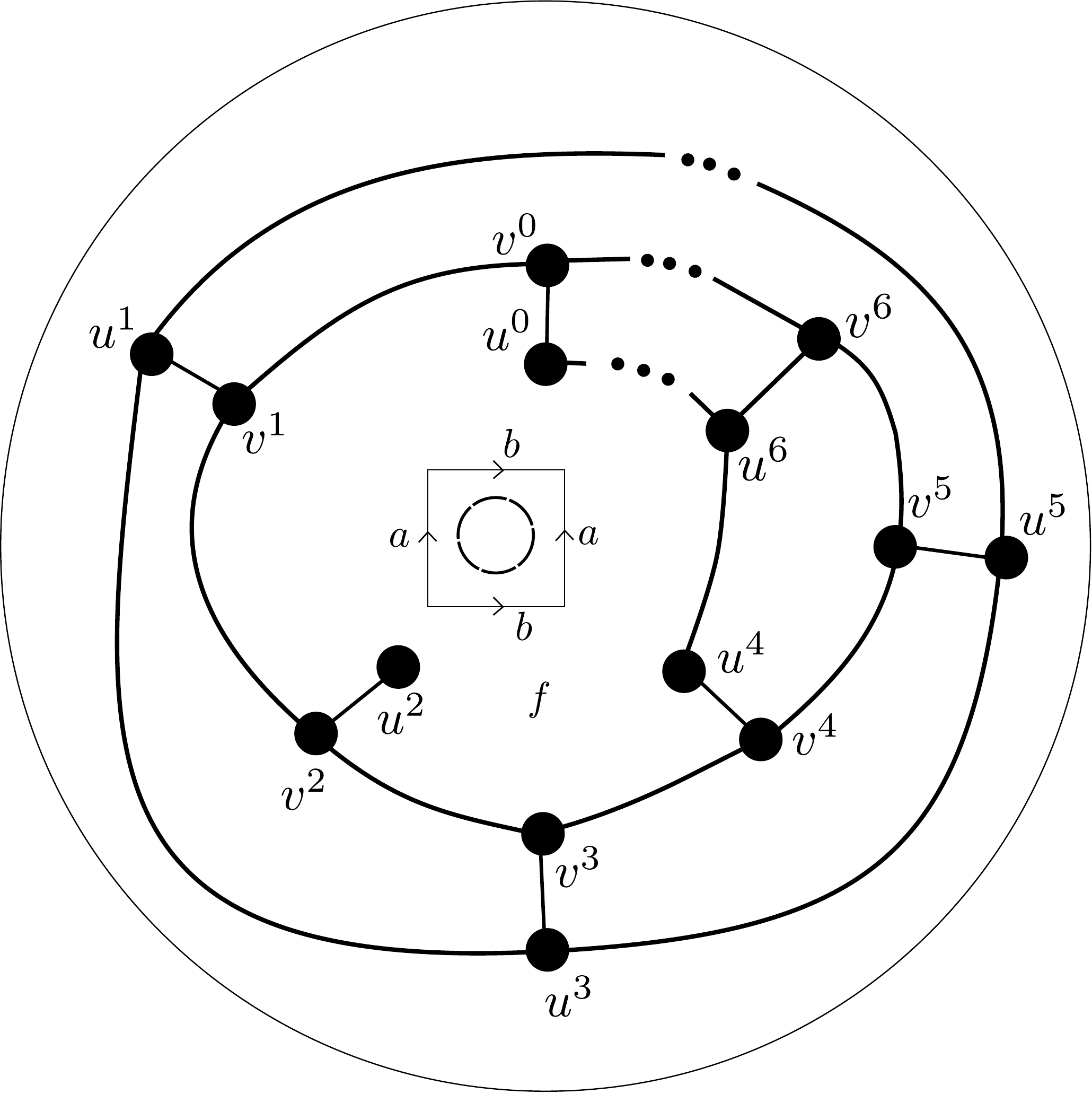}
\caption{A partial and noncelluar embedding of $GP(2p,p)$ in the torus.  The dashed circle in the middle of the square represents a boundary component of the interior of the square with edges identified.  The face $f$ is homeomorphic to a punctured torus.}\label{fig:GPStep2}
\end{center}
\end{figure}

Now, redraw the edges $(u^0,u^2)$ and $(u^2,u^4)$ such that they complete longitudinal and meridianal circles of the torus, as in Figure \ref{fig:GPStep3}.

\begin{figure}[H]
\begin{center}
\includegraphics[scale=.4]{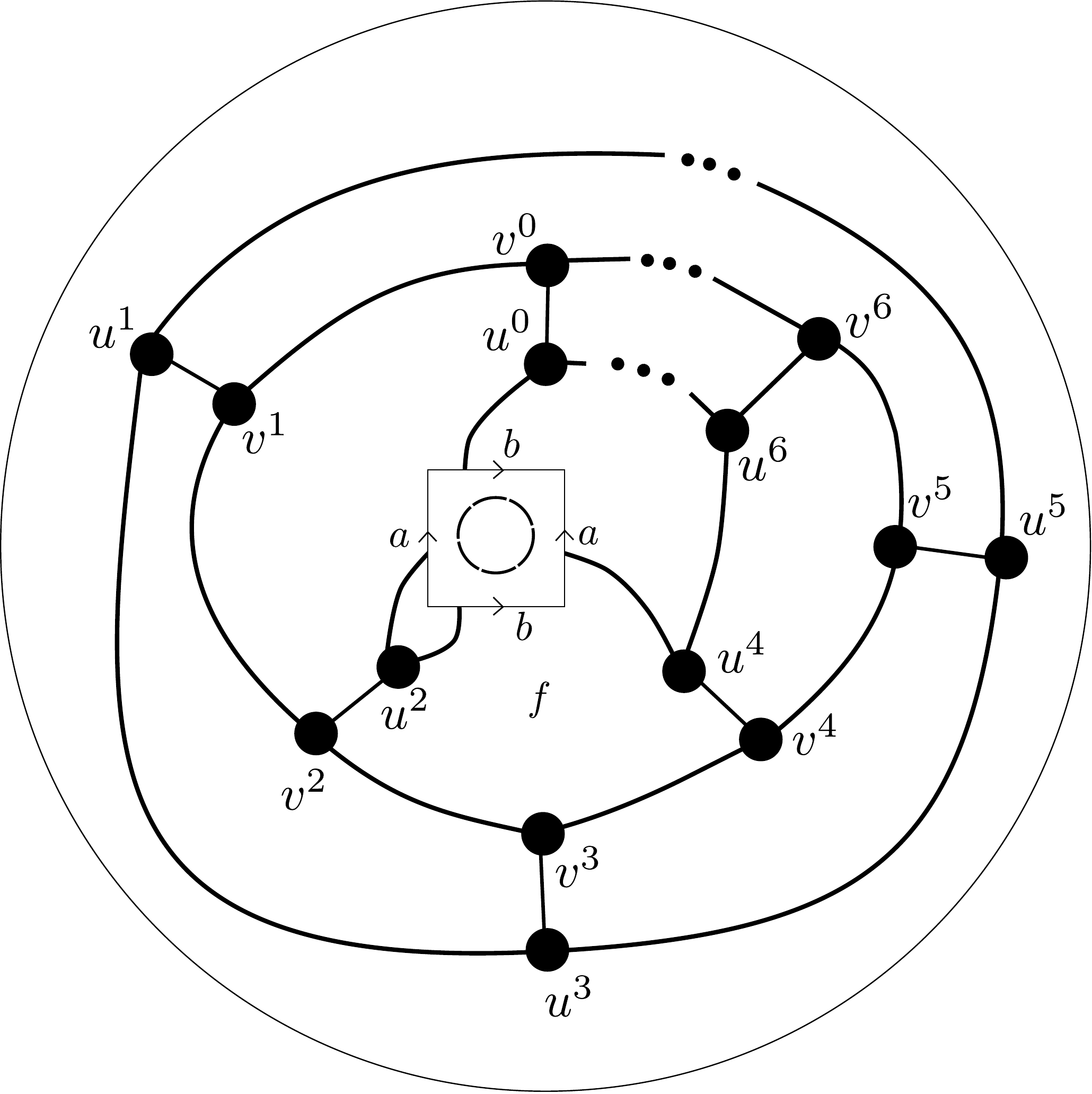}
\caption{$GP(2p,2)$ cellularly embedded in the torus.}\label{fig:GPStep3}
\end{center}
\end{figure}

It remains to check that the face $f$ shown in Figure \ref{fig:GPStep3} is homeomorphic to  a disc.  To verify this, note that in the case of each embedding of $GP(2p,2)$ in the torus constructed this way, a facial boundary of $f$ is always the same and can be described as a sequence of vertices, \[ u^0 v^0 v^1 v^2 u^2 u^4 u^6 \ldots u^0 u^2 v^2 v^3 v^4 u^4 v^4 u^2 u^0,\] where the vertices $\left \{u^i:\ i\ge 6\right \rbrace$ exist for $p>3$.  It is apparent that $f$ is homeomorphic to a disc.
\end{construction}

Construction \ref{constructon:GPTorus} completes the proof of Part \ref{theorem:ThePointPart2} of Theorem \ref{theorem:ThePoint}.\medskip

\noindent \textit{Proof of Part \ref{theorem:ThePointPart3}.}\medskip

\begin{figure}[H]
\begin{center}
\includegraphics[scale=.4]{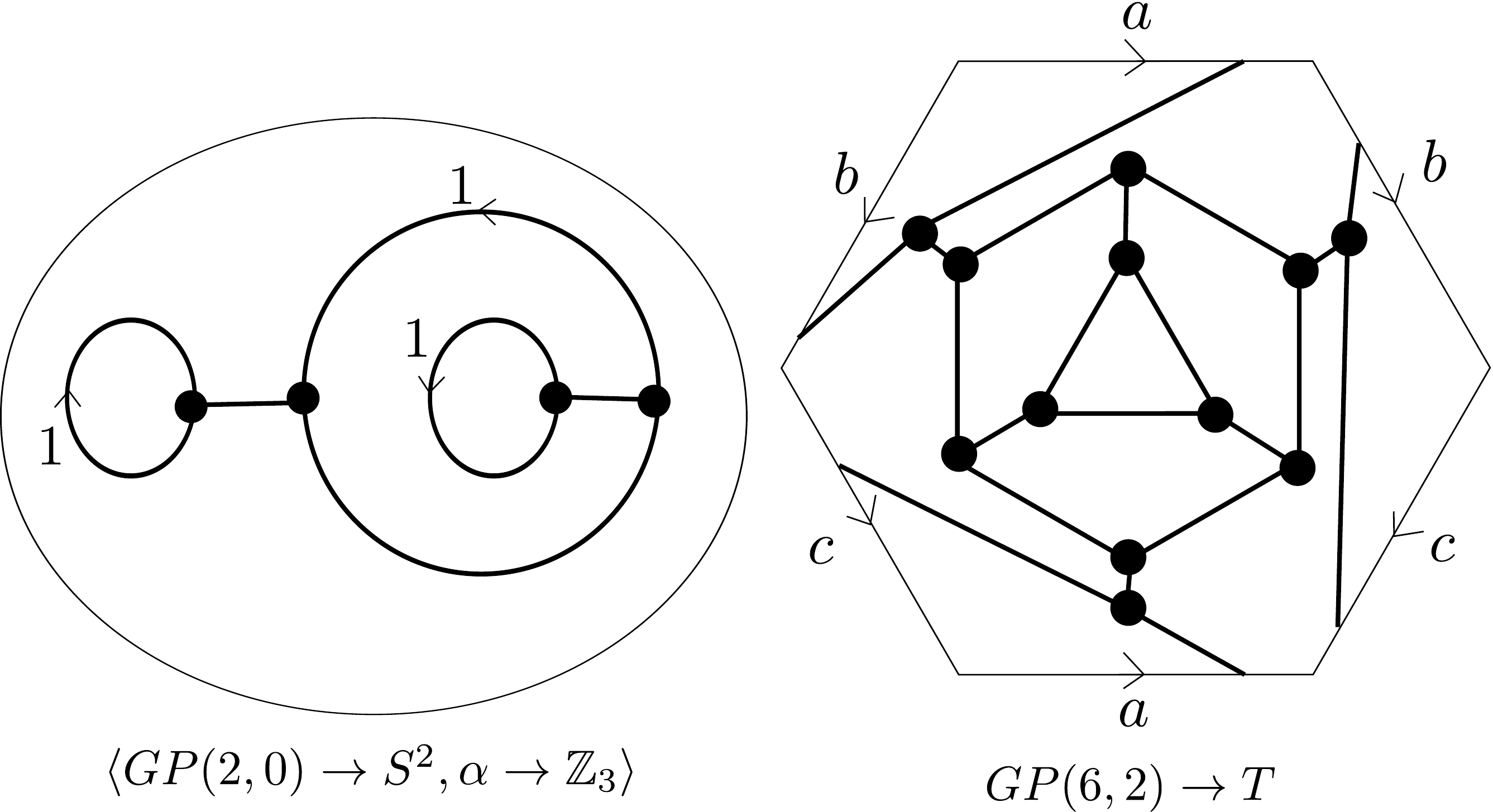}
\caption{An ordinary voltage graph and its derived embdding, which is $GP(6,2)\rightarrow T$.  The zero voltages are not shown, adn the darts on the derived graph are omitted.}\label{fig:GPBase5}
\end{center}
\end{figure}

We consider the graph $GP(2p,2)$ for which $p=3$. Figure \ref{fig:GPBase5} proves Part \ref{theorem:ThePointPart3} of Theorem \ref{theorem:ThePoint}.  It is easy to verify that the derived graph is indeed $GP(6,2)$.

\noindent \textit{Proof of Part \ref{theorem:ThePointPart4}.}\medskip

Now assume that $p$ is a prime greater than $5$.  If there were any automorphisms of $GP(2p,2)$ that map a vertex $v^j$ to a vertex $u^k$, then the $\mathbb{Z}_n$-action defined in Remark \ref{remark:ZnAction}, which is transitive on the vertices $u^i$ and the vertices $v^i$, respectively, can be combined with this automorphism to construct another automorphism of $GP(2p,2)$ that maps any vertex to any other vertex.  This would mean that $GP(2p,2)$ is vertex transitive, which contradicts a theorem of Frucht, Graver, and Watkins \cite[Theorem 1]{FGW}.  So, the only free-actions of groups on $GP(2p,2)$, for $p>5$ are those which have the vertices $v^i$ and $u^i$ in distinct orbits.  The vertices $v^i$ are the vertices of a circle of length $n$.  Since the automorphism group of a circle of length $n$ is is the dihedral group of order $2n$, and reflections of the circle have fixed points, it it follows that the only free actions of groups on $GP(2p,2)$ are the cyclic actions of the subgroups $\langle 1\rangle$, $\langle 2\rangle$, and $\langle p \rangle$ of $\mathbb{Z}_{2p}$, which act on $GP(2p,2)$ in the manner described in Remark \ref{remark:ZnAction}.  The quotient graphs of $GP(2p,2)$ modulo these actions are $GP(1,0)$, $GP(2,0)$, and $GP(p,2)$ respectively.  By Theorem \ref{theorem:VoltageGraphs}, there exist voltage assignments to these graphs such that the corresponding derived graphs are $GP(2p,2)$.  By \cite[Theorem 2.5.4]{GT}, we can require that $\alpha$ assigns voltage $0$ to the darts of a spanning tree of the ordinary voltage graph.  Thus, without loss of generality, the three possibilities for ordinary voltage graphs whose derived graphs are GP(2p,2) are: \begin{enumerate} 
\item{$\langle GP(1,0),\alpha \rightarrow \mathbb{Z}_{2p}\rangle$ such that $\alpha$ assigns voltage $1$ and $2$ to the positively directed edges on the loops $e_1$ $e_2$, respectively, and $0$ to both darts on the link, as in Figure \ref{fig:GPBase},}
\item{$\langle GP(2,0),\alpha,\rightarrow, \mathbb{Z}_p\rangle$, as in Figure \ref{fig:GPBase2},

\begin{figure}[H]
\begin{center}
\includegraphics[scale=.5]{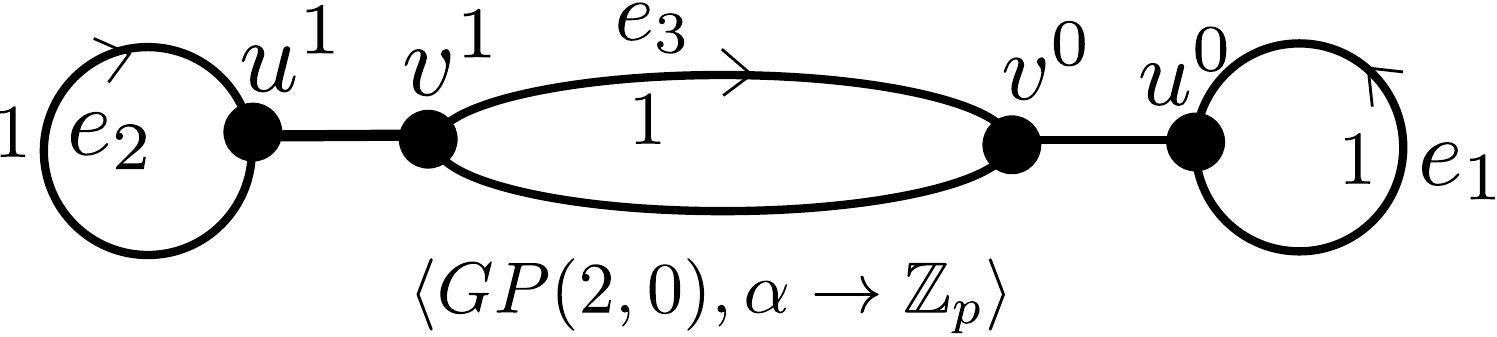}
\caption{An ordinary voltage graph whose derived graph is $GP(2p,2)$.  The zero voltages are not shown.}\label{fig:GPBase2}
\end{center}
\end{figure}
}

\item{$\langle GP(p,2),\alpha,\rightarrow, \mathbb{Z}_2\rangle$, such that $\alpha$ assigns $1$ to all darts on the edges joining the vertices $v^i$ and $0$ to all other darts.  An example of this appears in Figure \ref{fig:GPBase3}.

\begin{figure}[H]
\begin{center}
\includegraphics[scale=.5]{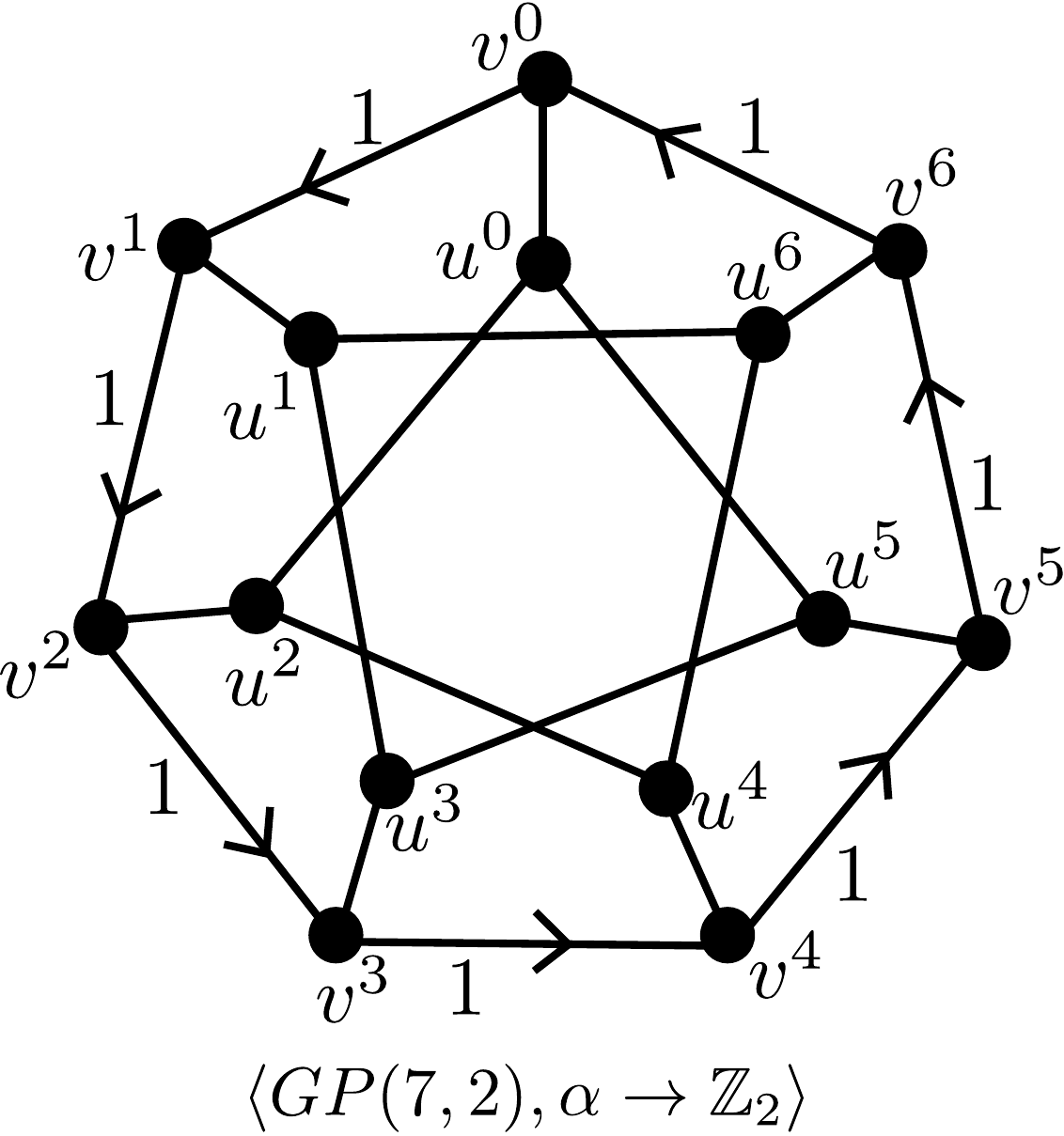}
\caption{An ordinary voltage graph whose derived graph is $GP(14,2)$.  The zero voltages are not shown.}\label{fig:GPBase3}
\end{center}
\end{figure}
}

\end{enumerate}

Lemma \ref{lemma:NoTorus} finishes the proof of Part \ref{theorem:ThePointPart4} of Theorem \ref{theorem:ThePoint}.

\begin{lemma}\label{lemma:NoTorus} For each prime $p>5$, the Generalized Petersen Graph $GP(2p,2)$ has no embedding in the torus as a derived embedding of an ordinary voltage graph embedding.\end{lemma}

\proof Considering the discussion in the proof of Part \ref{theorem:ThePointPart4} of Theorem \ref{theorem:ThePoint}, it follows from Theorems \ref{theorem:EulerCharConnectedSum}, \cite[Theorem 2.5.4]{GT}, and \ref{theorem:RHEquation} that we need only consider the three ordinary voltage graphs described in Figures \ref{fig:GPBase}, \ref{fig:GPBase2}, \ref{fig:GPBase3}, and all of their possible embeddings in the sphere, projective plane, torus, and Klein bottle.  We organize the proof by each ordinary voltage graph and then by the surface containing the base embedding.\medskip

\noindent\textbf{Case 1}: $GP(1,0)$\medskip

\noindent\textit{Case 1a}: Base embeddings of $GP(1,0)$ in the sphere\medskip

\begin{figure}[H]
\begin{center}
\includegraphics[scale=.5]{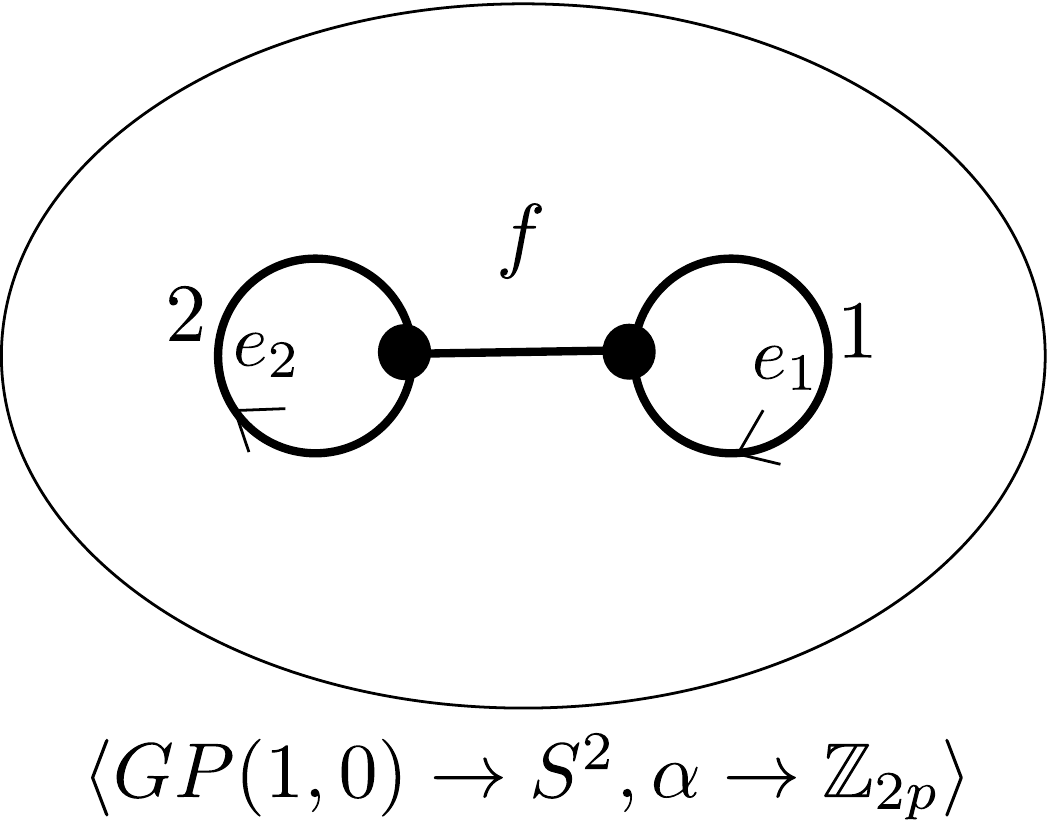}
\caption{An ordinary voltage graph embedding in the sphere.  The zero voltages are not shown.}\label{fig:GP11Sphere}
\end{center}
\end{figure}

Since $\mbox{dim}(Z(GP(1,0)))=2$ and the two loops induce circles that have no vertices in common, there is only one embedding with three faces: two faces are bounded by a loop, and the other by both loops and the link.  Figure \ref{fig:GP11Sphere} contains the only embedding of $GP(1,0)$ in the sphere though we will have to consider different orientations of the positive edges.

Given the embedding depicted in Figure \ref{fig:GP11Sphere}, note that the two faces bounded by the loops will always have facial boundary walks with net voltages of order $2p$ and $p$, respectively.  It follows that the deficiencies of the branch points contained in these two faces are $2p-1$ and $2p-2$, respectively.  The face $f$ has a facial boundary walk $W_f$ whose net voltage depends on whether the positive edges $e_1$ and $e_2$ (or both of their opposites) appear in $W_f$.  If this is true, which is the case depicted in Figure \ref{fig:GP11Sphere}, then $\omega(W_f)$ is either $3$ or $-3$, neither of which divides $2p$ for $p>5$.  It follows by Theorem \ref{theorem:RHEquation} that \[ \chi(S^2)^\alpha = 4p - ((2p-1)+(2p-1)+(2p-2))=4-2p<0,\] which precludes $(S^2)^\alpha$ from being the torus.  The only other possibility is that $W_f$ contains $e_1$ and $e_2^{-}$ (or vice versa).  In this case it is apparent that \[\langle \omega(W_f)\rangle = \mathbb{Z}_{2p},\] so the deficiency of the branch point contained in $f$ is $2p-1$.  Just as before, $\chi((S^2)^\alpha)<0$, which precludes $(S^2)^\alpha$ from being the torus.\medskip

\noindent\textit{Case 1b}: Base embeddings of $GP(1,0)$ in the projective plane\medskip

Consider the embedding of $GP(1,0)$ appearing in the base embedding in Figure \ref{fig:SphereDerived}.  First note that at least one of $e_1$ and $e_2$ must represent the nontrivial homology class of $P^2$.  Since $e_1$ and $e_2$ induce circles in $GP(1,0)$ have no vertices in common, Lemma \ref{lemma:P2CirclesCross} implies that only one of the loops can be homologically nontrivial, and so, the other loop must bound a face.  Therefore, there is no other cellular embedding of $GP(1,0)$ in $P^2$ to consider.

We now determine the deficiency of any branch points.  By Lemma \ref{lemma:OrientationReversingVoltage}, it cannot be true that $[e_2]\neq [0]$ and $(P^2)^\alpha$ is orientable.  It follows that the face bounded by a loop must contain a branch point of deficiency $2p-2$.  The other face has a facial boundary walk whose net voltage depends on whether the positive edge $e_1$ and $e_2$ (or both of their opposites) appear in a facial boundary walk $W_f$.  If they do, then $\omega(W_f)$ generates the subgroup $\langle 4\rangle$ of $\mathbb{Z}_{2p}$.  In $\mathbb{Z}_{2p}$, $\langle 2\rangle = \langle 4\rangle$ since $2$ and $4$ have the same greatest common divisor with $2p$.  It follows that both faces contain branch points of deficiency $2p-2$, and by Theorem \ref{theorem:RHEquation}, we have \[\chi((P^2)^\alpha) = 2p - ((2p-2)+(2p-2)) = 4 - 2p,\] which precludes $(P^2)^\alpha$ from being the torus.  The other case where $e_1$ and $e_2^-$, or their respective opposites appear in $W_f$, was treated in Figure \ref{fig:SphereDerived} and Remark \ref{remark:SphereDerived} and always produces the sphere as the derived surface.\medskip

\noindent\textit{Case 1c:} Base embeddings of $GP(1,0)$ in the Klein bottle\medskip

Recall that $\beta_1(\mbox{\textit{KB}})=2$ and note that $\mbox{dim}(Z(GP(1,0)))=2$.  So, $e_1$ and $e_2$ must be homologically nontrivial and homologically independent.  We will show that $e_2$ induces an orientation-reversing circle, which, by Lemma \ref{lemma:OrientationReversingVoltage}, will preclude $\mbox{\textit{KB}}^\alpha$ from being the torus.  Since the Klein bottle is nonorientable, there must be at least one orientation-reversing circle; assume that it is induced by $e_1$.  Since $e_2$ must be homologically nontrivial, if it does not induce an orientation-reversing circle, it must induce a nonseparating orientation-preserving circle.  If $e_2$ induces a nonseparating orientation-preserving circle, Lemma \ref{lemma:KBCirclesCross} implies that the circles induced by $e_1$ and $e_2$ must transversely cross eachother, which is impossible since they have no vertices in common.  It follows that $e_2$ must induce an orientation-reversing circle for any cellular embedding of $GP(1,0)$ in the Klein bottle.\medskip

\noindent\textit{Case 1d}: Base embeddings of $GP(1,0)$ in the torus\medskip

Recall that $\beta_1(T)=2$ and note again that $\mbox{dim}(Z(GP(1,0)))=2$.  So $e_1$ and $e_2$ must be homologically nontrivial and homologically independent.  By Lemma \ref{lemma:TorusCirclesCross}, $\langle e_1,e_2\rangle =1$, which is impossible since the circles induced by $e_1$ and $e_2$ have no vertices in common.  Therefore, there are no cellular embeddings of $GP(1,0)$ in the torus.\medskip

\noindent\textbf{Case 2}: $GP(2,0)$\medskip

\noindent\textit{Case 2a}: Base embeddings of $GP(2,0)$ in the sphere\medskip

\begin{figure}[H]
\begin{center}
\includegraphics[scale=.4]{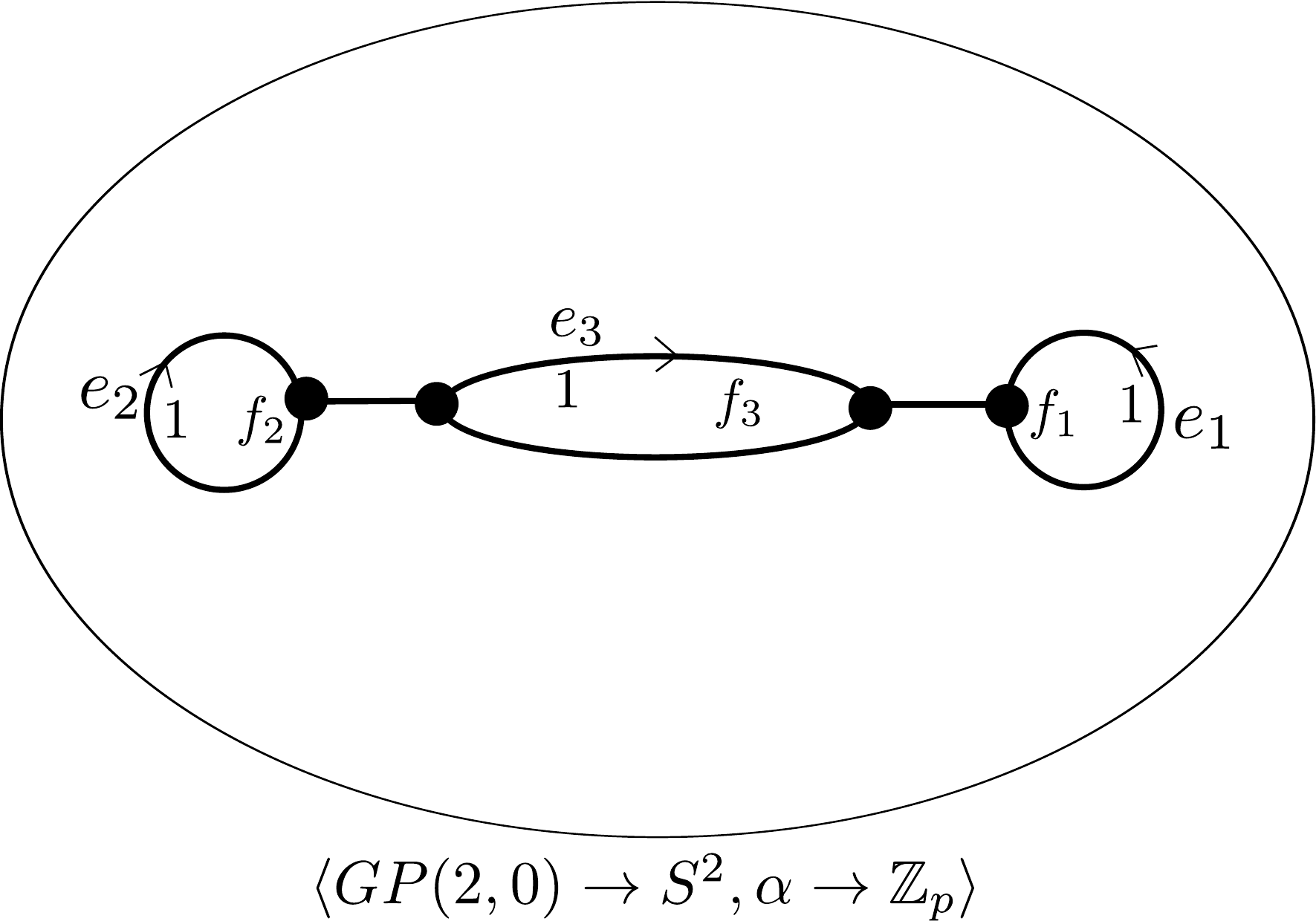}	
\caption{An ordinary voltage graph embedding featuring one of the two possible embeddings of $GP(2,0)$ in the sphere.  The zero voltages are not shown.}\label{fig:GPBase2a}
\end{center}
\end{figure}

\begin{figure}[H]
\begin{center}
\includegraphics[scale=.4]{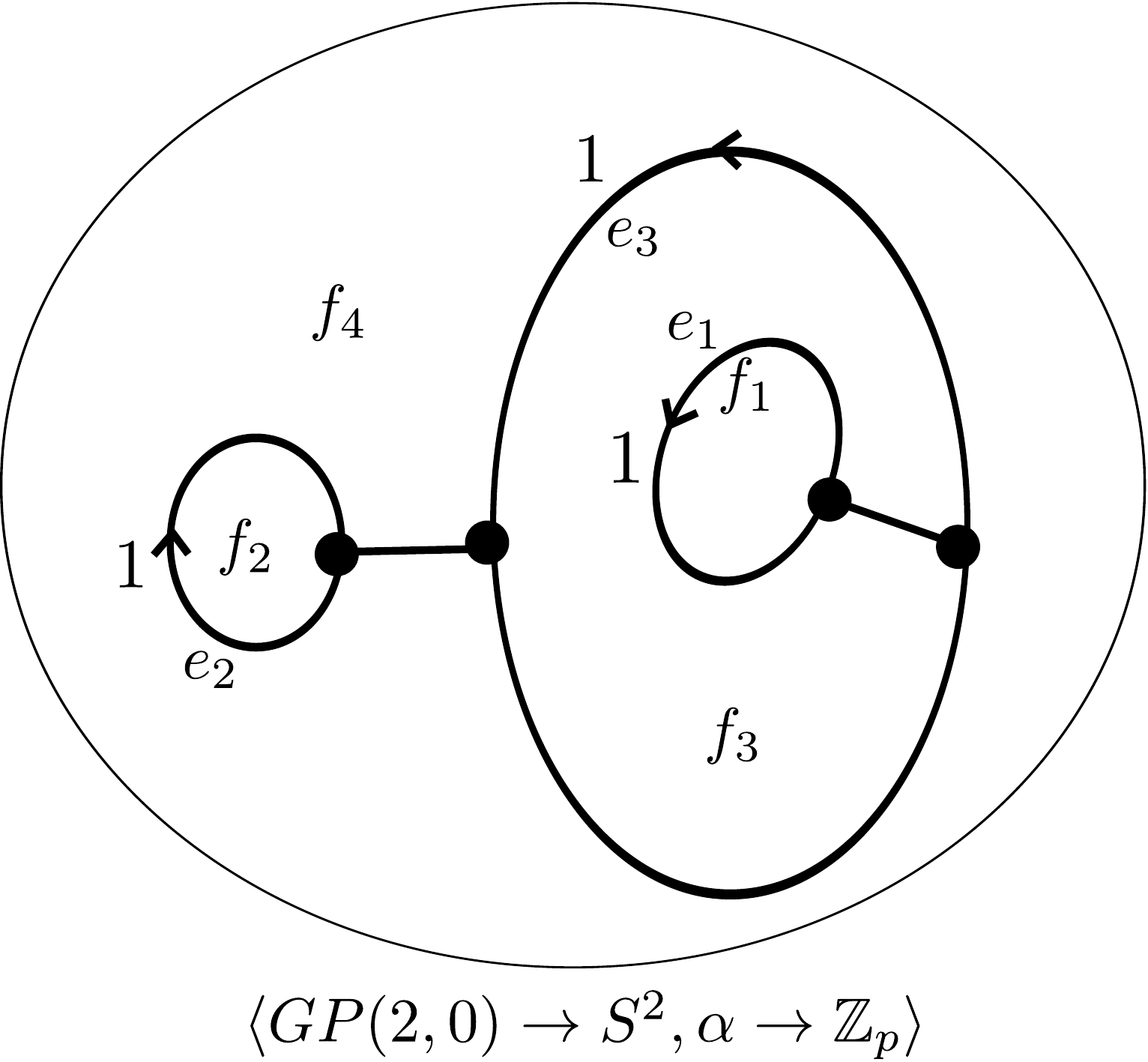}
\caption{An ordinary voltage graph embedding featuring one of two possible embeddings of $GP(2,0)$ in the sphere.}\label{fig:GPBase2b}
\end{center}
\end{figure}

Since any two of the three circles of $GP(2,0)$ have no vertices in common, there can be no crossings of the circles induced by any of the nonzero elements of $Z(GP(2,0)$.  The Jordan Curve Theorem guarantees that all circles embedded in the sphere are separating, and so it follows that there are only two embeddings of $GP(2,0)$ in the sphere: one of the embeddings features both loops contained in a single component of the complement of the circle of length 2, and the other has one loop contained in each of those components.  These cases are drawn in the base embeddings in Figures \ref{fig:GPBase2a} and \ref{fig:GPBase2b}, and we first consider the former case.

No matter how the positive directions on $e_1$ and $e_2$ are drawn in this embedding, each of the three faces $f_1$, $f_2$ and $f_3$	contain branch points of deficiency $p-1$. By Theorem \ref{theorem:RHEquation}, it follows that \[\chi((S^2)^\alpha)\le 2p - ((p-1)+(p-1)+(p-1))=3-p.\]  Since $p>5$, this precludes $(S^2)^\alpha$ from being the torus.

Now consider Figure \ref{fig:GPBase2b}.

No matter how the positive edges $e_1$ and $e_2$ are drawn in this embedding, each of the faces $f_1$ and $f_2$ will have facial boundary walks whose net voltage is $1$, which generates $\mathbb{Z}_p$.  In the depicted case, the facial boundary walks of $f_3$ and $f_4$ will have net voltage $0$, and so $f_3$ and $f_4$ have no branch points.  It follows that $(S^2)^\alpha=S^2$ since Theorem \ref{theorem:RHEquation} implies that \[\chi((S^2)^\alpha) = 2p =((p-1) + (p-1)) =2.\]

However, if any of the positive directed edges $e_1$, $e_2$, and $e_3$ were to be drawn in such a way that at least one of the faces $f_3$ or $f_4$ would have a nonzero voltage, one of them would have a facial boundary walk with net voltage $2$, which generates $\mathbb{Z}_p$, then Theorem \ref{theorem:RHEquation} implies \[\chi((S^2)^\alpha)\le 2p - ((p-1)+(p-1)+(p-1))=3-p.\]  Since $ p > 5$, this precludes $(S^2)^\alpha$ from being the torus.\medskip

\noindent\textit{Case2b}: Base embeddings of $GP(2,0)$ in the projective plane and the Klein bottle\medskip

In either of these surfaces, there must be an orientation-reversing circle.  There are only three circles in $GP(2,0)$, and each of them is traversed by an Eulerian walk having net voltage $1$, which has odd order in $\mathbb{Z}_p$.  It follows by Lemma \ref{lemma:OrientationReversingVoltage} that any base embedding of $GP(2,0)$ in the projective plane or the Klein bottle with the voltage assignment described in Figure \ref{fig:GPBase2} results in a nonorientable derived surface.\medskip

\noindent\textit{Case2c}: Base embeddings of $GP(2,0)$ in the torus\medskip

Since each of the three circles of $GP(2,0)$ have no vertices in common, Lemma \ref{lemma:TorusCirclesCross} excludes the possibility of there being any base embeddings of $GP(2,0)$ in the torus.\medskip

\noindent\textbf{Case 3}: $GP(p,2)$\medskip

Let $C_v,\ C_u\in Z(GP(p,2))$ denote the 1-chains inducing the circles whose vertex sets are $\left\{ v_i\right \rbrace$ and $\left \{u_i\right \rbrace$, respectively.  Since $2$ does not divide $p$, there is a single circle connecting the vertices $u_i$.  It follows that $GP(p,2)\setminus (GP(p,2)\mycolon C_v)$ is connected.  For future use we define some fundamental cycles with respect to a spanning tree of $GP(p,2)$.  We chose the spanning tree \[GP(p,2)\mycolon ( \left \{ (v_i,v_{i+1}):\ i\in \mathbb{Z}_{p-1} \right \rbrace \bigcup \left \{(v_i,u_i):\ i \in \mathbb{Z}_p\right \rbrace,\] which leaves us with the fundamental edges \[ \left \{(v_{p-1},v_0)\right \rbrace \cup \left \{ (u_i,u_{i+2})\ :\ i\in \mathbb{Z}_p\right \rbrace,\]  and, for each $i\in \mathbb{Z}_p$ we define the fundamental cycle \[C_i = (u_i,u_{i+2})+(v_{i+2},u_{i+2})+(v_{i+1},v_{i+2})+(v_i,v_{i+1})+(v_i,u_i).\]  The fundamental cycles $C_i$ and $C_v$ together form the set of fundamental cycles with respect to our spanning tree.  Moreover, for each $i\in \mathbb{Z}_p$, an Eulerian walk $W_i$ of $GP(p,2)\mycolon C_i$ satisfies $\omega(W_i)=0$ since $A=\mathbb{Z}_2$ and we are assigning a voltage assignment that assigns voltage $1$ to exactly two darts in $W_i$.  It follows by Part \ref{theorem:BigCosetTheoremPartThree} of Theorem \ref{theorem:BigCosetTheorem} that there are exactly two circles forming the fiber over each circle $GP(p,2):C_i$.  Making shorthand of our previous notation, we will let $C_i^0$ and $C_i^1$ denote the 1-chains forming the fiber over $C_i$, according to whether their induced circles contain $v_i^0$ or $v_i^1$, respectively.\medskip

\noindent\textit{Case 3a}: $GP(p,2)$ in the sphere\medskip

We will show that $GP(p,2)$ has a $K_{3,3}$ minor, which, by Kuratowski's theorem, precludes $GP(p,2)$ from being planar.  Let $X=\left \{v_0,v_2,u_1\right \rbrace$ and $Y=\left \{v_1,v_3,u_2\right \rbrace$. We identify $9$ edge sets inducing 9 $X$-$Y$ independent paths.

\begin{enumerate}

\item{Edge sets inducing paths joining $v_1$ to all vertices in $X$ are

		\begin{enumerate}
		
		\item{$\left \{ (v_1,v_0)\right \rbrace$ inducing a path joining $v_1$ and $v_0$, and}
		
		\item{$\left \{(v_1,u_1)\right \rbrace$ inducing a path joining $v_1$ and $u_1$, and}
		
		\item{$\left \{(v_1,v_2)\right \rbrace$ inducing a path joining $v_1$ and $v_2$.}
		
		\end{enumerate}
		}
		
\item{Edge sets inducing paths joining $v_3$ to all vertices in $X$ are:

		\begin{enumerate}

		\item{$\left \{(v_3,v_4),\ (v_4,v_5),\ldots, (v_{p-1},v_0)\right \rbrace$ inducing a path joining $v_3$ and $v_0$.}

		\item{$\left \{(v_3,u_3),\ (u_3,u_1)\right \rbrace$ inducing a path joining $v_3$ and $u_1$, and}
		
		\item{$\left \{(v_3,v_2)\right \rbrace$ inducing a path joining $v_3$ and $v_2$.}

		\end{enumerate}
		}
\item{Edge sets inducing paths joining $u_2$ to all vertices in $X$ are:

		\begin{enumerate}		

		\item{$\left \{(u_2,u_4), (u_4,u_6), \ldots, (u_{p-1},u_1) \right \rbrace$ inducing a path joining $u_2$ and $u_1$, and}

		\item{$\left \{(u_2, u_0),\ (u_0,v_0)\right \rbrace$ inducing a path joining $u_2$ and $v_0$, and}
		
		\item{$\left \{(u_2,v_2) \right \rbrace$ inducing a path joining $u_2$ and $v_2$.}

		\end{enumerate}
		
		}

\end{enumerate}It follows that $K_{3,3}$ is a minor of $GP(p,2)$.  Therefore, there is no ordinary voltage graph embedding of $GP(p,2)$ in the sphere.\medskip

\textit{Case 3b}: $GP(p,2)$ in the projective plane\medskip

In this case, since each dart on an edge appearing in $C_u$ has voltage $0$, it follows by Lemma \ref{lemma:OrientationReversingVoltage}, that for $(P^2)^\alpha$ to be the torus, $C_u$ must be homologically trivial.  Similarly, each $C_i$ must also be homologically trivial.  Lemma \ref{lemma:P2CirclesCross} implies that $C_u$ and each $C_i$ induces a contractible circle.  We organize the remainder of these cases according to whether $C_v$ is homologically trivial.\medskip

\textit{Subcase 1}: $C_v$ is homologically nontrivial\medskip

\begin{figure}[H]
\begin{center}
\includegraphics[scale=.45]{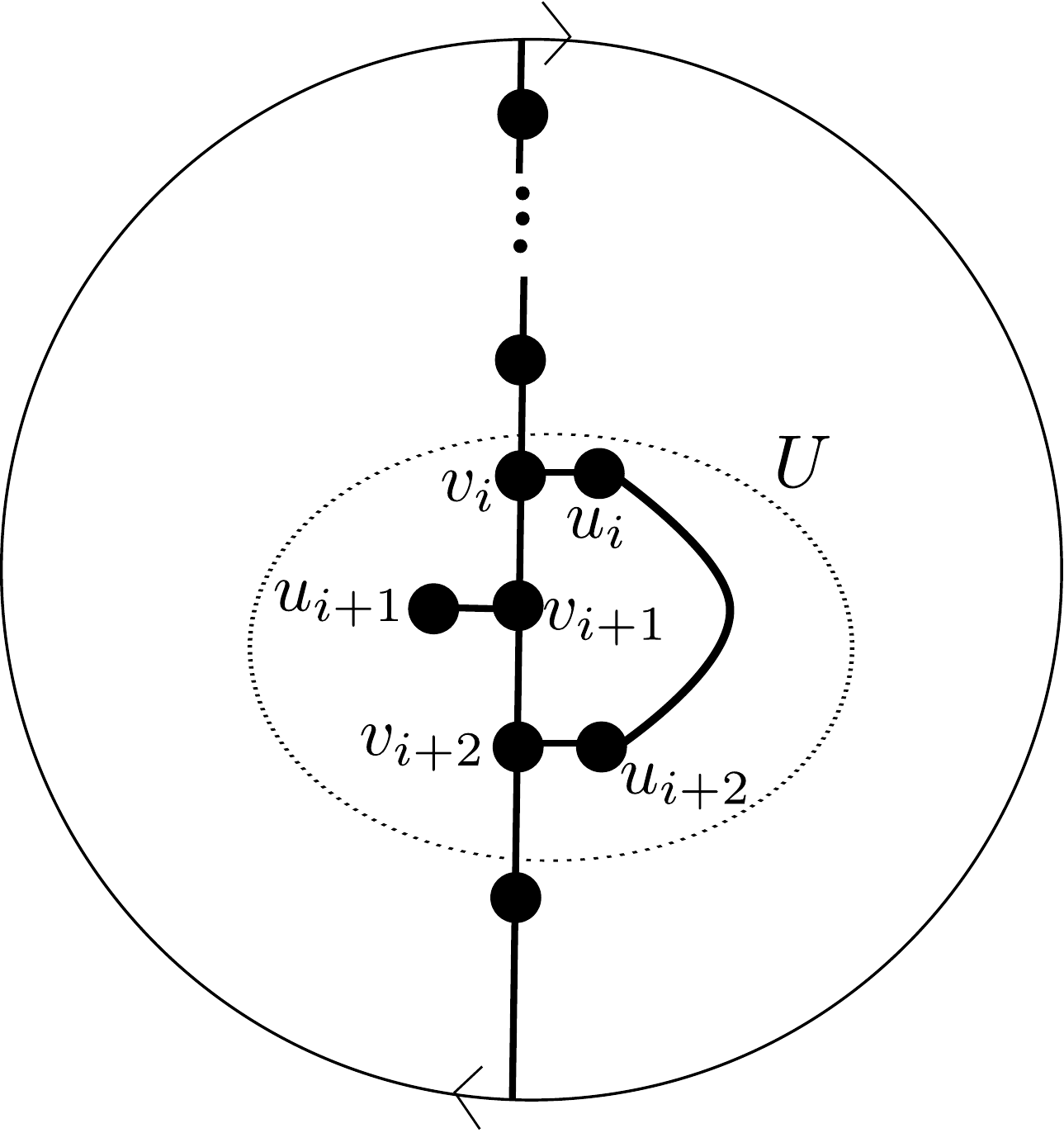}
\caption{A partial embedding of $GP(p,2)$ in $P^2$.}\label{fig:GPp2P2Base1}
\end{center}
\end{figure}

If the circle $GP(p,2)\mycolon C_i$ is orientation preserving, then $GP(p,2)\mycolon C_i$ is contractible in $P^2$, and we must have the edge $(v_{i+1},u_{i+1})$ appearing as it does in the disc $U$ as drawn in Figure \ref{fig:GPp2P2Base1}. 

\begin{figure}[H]
\begin{center}
\includegraphics[scale=.45]{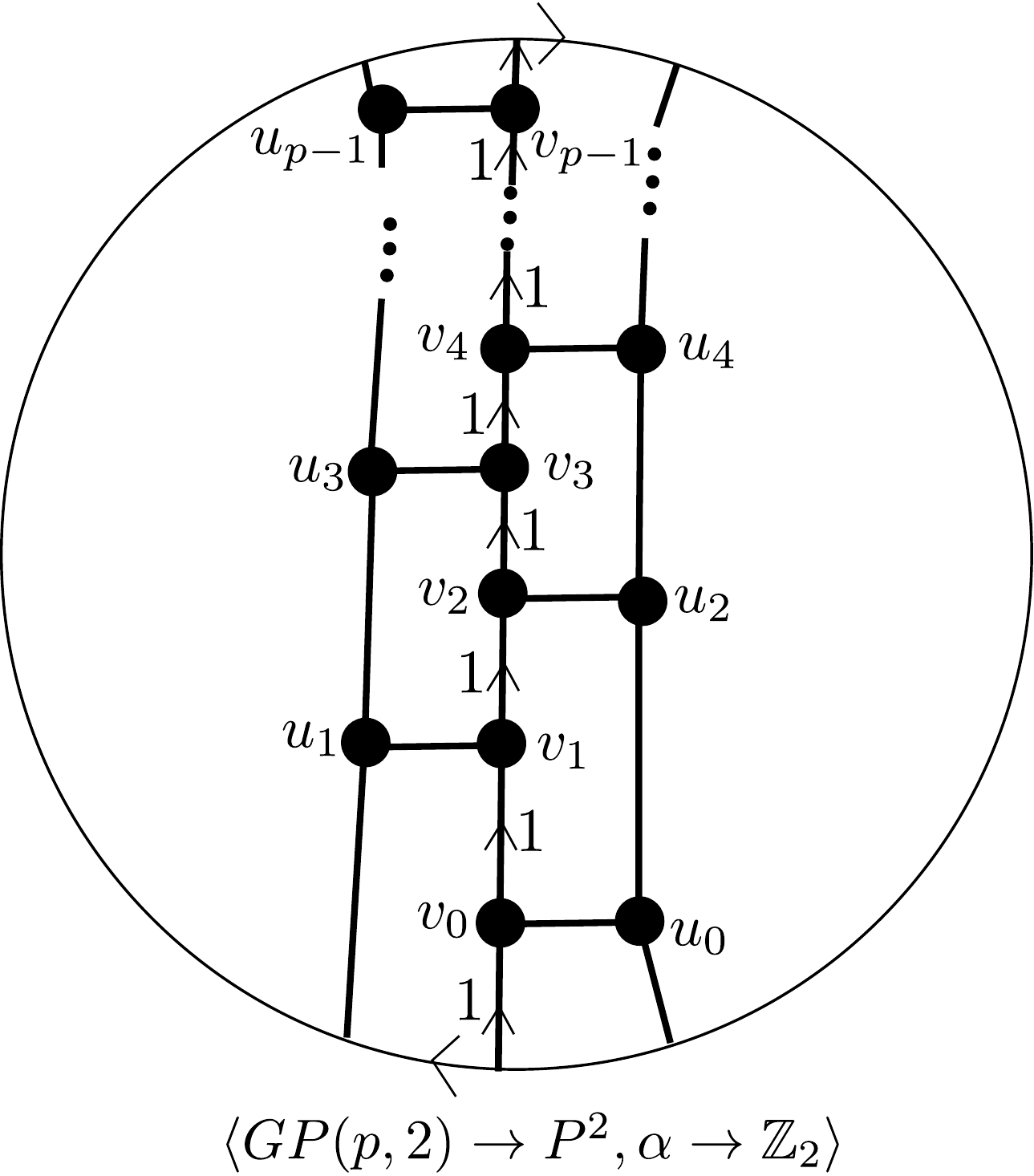}
\caption{An ordinary voltage graph embedding in the projective plane.  The zero voltages are not shown.}\label{fig:GPp2P2Base3}
\end{center}
\end{figure}

Since each $C_i$ induces a contractible circle, an embedding of this type must be the embedding depicted in Figure \ref{fig:GPp2P2Base3}.

 Since each face of the embedding in Figure \ref{fig:GPp2P2Base3} is bounded by an even number of edges contained in $C_v$, it follows that there are no branch points.  It follows from Theorem \ref{theorem:RHEquation} that \[\chi((P^2)^\alpha) = 2.\] And so, $(P^2)^\alpha$ is not the torus.\medskip

\textit{Subcase 2}: $C_v$ is homologically trivial\medskip

\begin{figure}[H]
\begin{center}
\includegraphics[scale=.45]{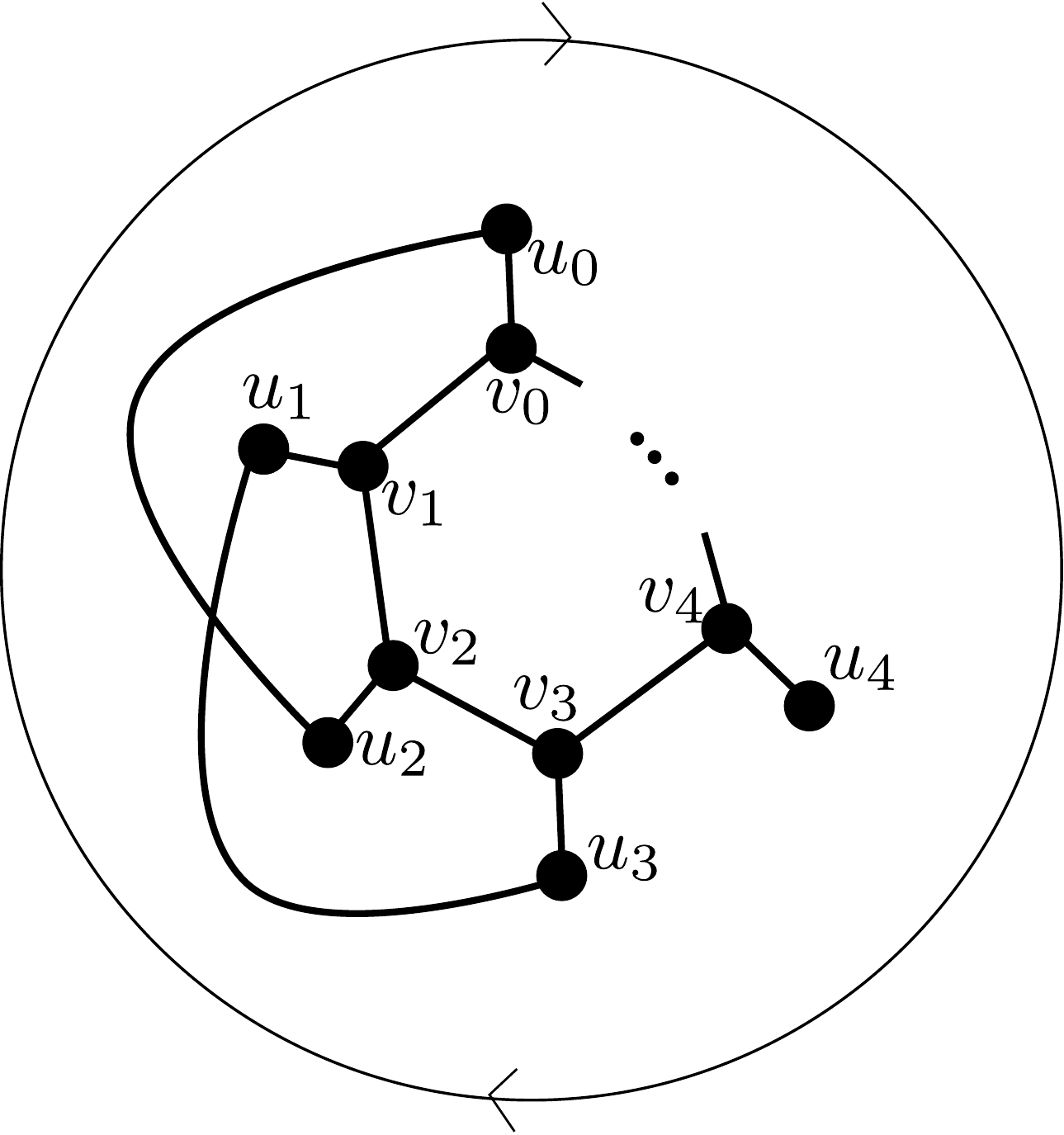}
\end{center}
\caption{A partial drawing of $GP(p,2)$ in the projective plane depicting a forced immersion of $GP(p,2)$.}\label{fig:GPp2P2Base4} 
\end{figure}

In this case, Lemma \ref{lemma:P2Decomposition} implies that $GP(p,2)\mycolon C_v$ is a contractible circle.  Since $GP(p,2)\mycolon C_v$ is connected and $GP(p,2)$ was determined to be nonplanar, $GP(p,2)\setminus E(C_v)$ must not be contained in the disc bounded by $GP(p,2)\mycolon C_v$.  In Figure \ref{fig:GPp2P2Base4}, a subgraph of $GP(p,2)$ is drawn in a manner that reflects the above constraints.  Since $C_0$ is assumed to be homologically trivial, Lemma \ref{lemma:P2Decomposition} implies that $GP(p,2)\mycolon C_0$ bounds a disc containing the vertex $u_1$.  If follows that if $C_0$ is homologically trivial then $GP(p,2)$ cannot be embedded in the projective plane.

It follows that $C_0$ must be homologically nontrivial in this case.  By Lemma \ref{lemma:OrientationReversingVoltage}, $(P^2)^\alpha$ is nonorientable, which precludes $(P^2)^\alpha$ from being the torus.\medskip

\noindent\textit{Case 3c}: $GP(p,2)$ in the Klein bottle\medskip

\noindent\textit{Subcase1}: $C_v$ is homologically trivial\medskip

If $GP(p,2)\mycolon C_v$ is a separating curve, then the fact that $GP(p,2)\setminus E(C_v)$ is connected means that that $GP(p,2)\setminus GP(p,2)\mycolon C_v$ is contained within one of the two regions of $\mbox{\textit{KB}}\setminus GP(p,2)\mycolon C_v$.  By Lemma \ref{lemma:KBDecomposition}, $C_v$ induces either a contractible circle forming the boundary of the punctured Klein bottle and the punctured sphere in the connected sum $\mbox{\textit{KB}}=\mbox{\textit{KB}}\#S^2$, or it induces a circle forming the boundary of the two projective planes forming the connected sum $\mbox{\textit{KB}}=P^2\#P^2$.

If $GP(p,2)\mycolon C_v$ is contractible, the fact that the $GP(p,2)$ is nonplanar means that $GP(p,2)\setminus E(C_v)$ is contained within the closure of the punctured Klein bottle, in which case $GP(p,2)\mycolon C_v$ bounds a face $f_v$.  Since $p$ is odd, $A=\mathbb{Z}_2$, and each of the darts (positive or negative) on the edges of the $GP(p,2)\mycolon C_v$ is assigned voltage $1$, it follows that a facial-boundary walk $W_f$ of $f$ satisfies $\omega(W_f)=1$.  Since the voltage group $A=\mathbb{Z}_2$, it follows that $f$ contains a branch point $y$ of deficiency $1$.  Therefore, by Theorem \ref{theorem:RHEquation}, \[ \chi (\mbox{\textit{KB}}^\alpha)\le 2\cdot 0 -1=-1, \] which precludes $\mbox{\textit{KB}}^\alpha$ from being the torus.

If $GP(p,2)\mycolon C_v$ is the boundary of the two punctured projective planes forming the connected sum $\mbox{\textit{KB}}=P^2\#P^2$, then the fact that $GP(p,2)\setminus E(C_v)$ is connected means that $GP(p,2)\setminus GP(p,2)\mycolon C_v$ must be contained within one of the two punctured projective planes.  It follows that the other punctured projective plane is a face of the embedding, which means that in this instance, the embedding is not cellular, a contradiction since ordinary voltage graph embeddings are by definition cellular embeddings.\medskip

\noindent\textit{Subcase 2}: $C_v$ is homologically nontrivial\medskip

Either $GP(p,2)\mycolon C_v$ is an orientation-preserving circle or $GP(p,2)\mycolon C_v$ is an orientation-reversing circle.

Suppose that $GP(p,2)\mycolon C_v$ is orientation-preserving.  Since $GP(p,2)$ is assumed to be cellularly embedded in a nonorientable surface, there must be at least one orientaton-reversing circle induced by one of the $C_i$, say $C_0$.  For aforementioned reasons, it follows from Lemma \ref{lemma:OrientationReversingVoltage} that $\mbox{\textit{KB}}^\alpha$ is nonorientable, so it is not the torus.

\begin{figure}[H]
\begin{center}
\includegraphics[scale=.45]{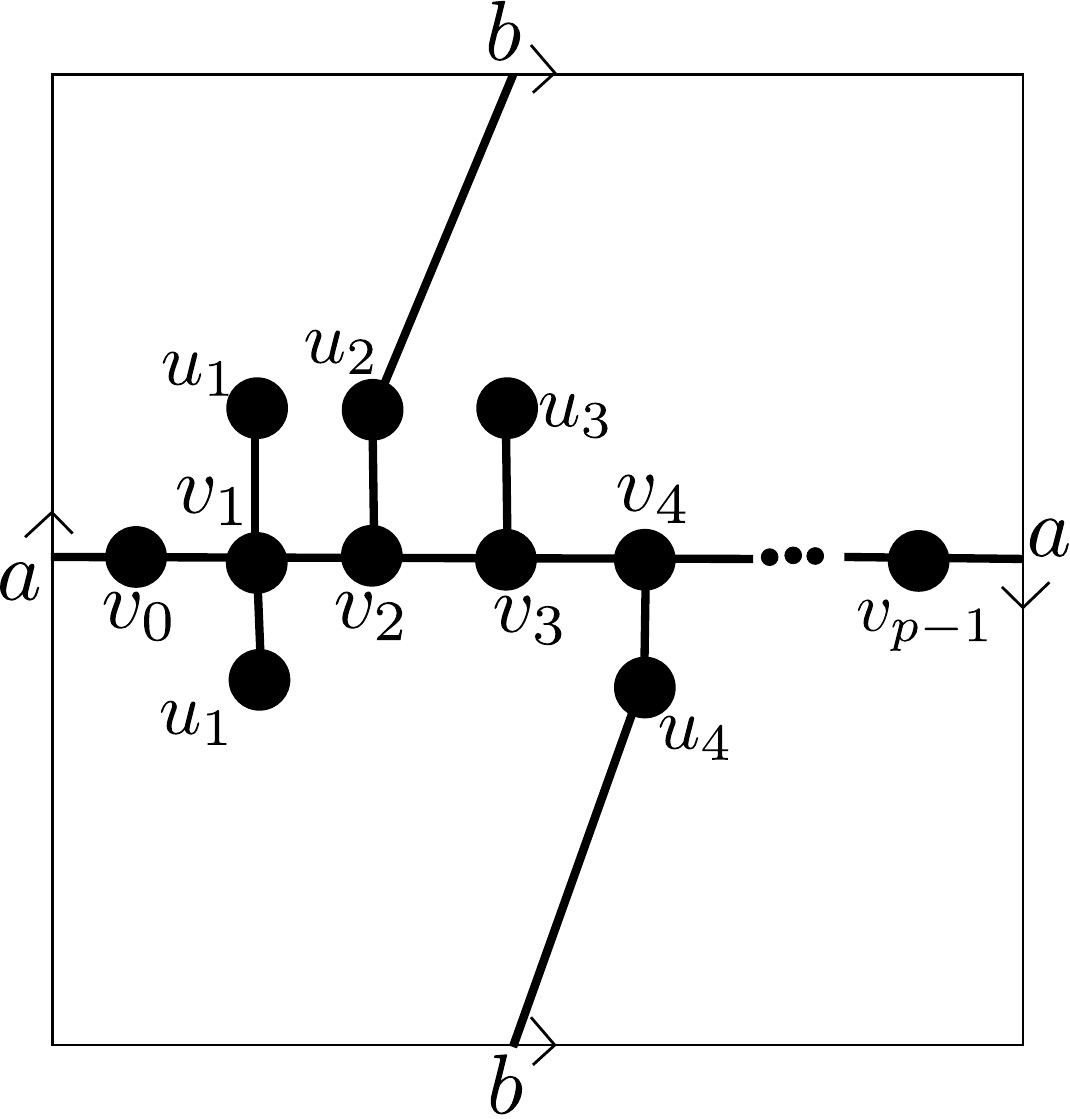}
\caption{A partial embedding of $GP(p,2)$ in the Klein bottle, with two possible choices for how to draw $(v_1,u_1)$.}\label{fig:GPp2KBBase1}
\end{center}
\end{figure}

Now, suppose that $GP(p,2)\mycolon C_v$ is an orientation-reversing circle.  Since $\beta_1(\mbox{\textit{KB}})=2$ and every element of $Z(G)$ is generated by $C_v$ and the $C_i$, it follows that at least one of the $C_i$, say $C_2$ is homologically nontrivial.  If $GP(p,2)\mycolon C_2$ is orientation-reversing then it follows by Lemma \ref{lemma:OrientationReversingVoltage} that $\mbox{\textit{KB}}^\alpha$ is nonorientable, in which case $\mbox{\textit{KB}}^\alpha$ is not the torus.  Therefore $GP(p,2)\mycolon C_2$ must be orientation preserving.  By Lemma \ref{lemma:KBCirclesCross}, $GP(p,2)\mycolon C_2$ may be drawn as it is in Figure \ref{fig:GPp2KBBase1}, transversely crossing $GP(p,2)\mycolon C_v$.  Note that without loss of generality, $(v_3,u_3)$ may be drawn as it is in Figure \ref{fig:GPp2KBBase1}, above $GP(p,2)\mycolon C_v$.  Also note that the two possible placements for the edge $(u_1,v_1)$ are also shown.

No matter how we draw the edge $(u_1,v_1)$, the fundamental cycles $C_1$ and $C_2$ must satisfy $\langle C_1,C_2\rangle=1$.  Since $[C_v]$ and $[C_2]$ span the $\mathbb{Z}_2$-vector space $H_1(\mbox{\textit{KB}})$, $C_1$ is homologous to $C_v$, $C_2$ or $C_v+C_2$.  Since $\langle C_v,C_2\rangle=1$, it follows that $C_1$ and $C_2$ are not homologous since $GP(p,2)\mycolon C_v$ is an orientation-preserving circle, which implies that $\langle C_2,C_2\rangle=\langle C_1,C_2\rangle=0$.  If $C_1$ is homologous to $C_v$, then Remark \ref{remark:TransverseCrossings} implies that $GP(p,2)\mycolon C_v$ is an orientation-reversing circle, and Lemma \ref{lemma:OrientationReversingVoltage} implies that $\mbox{\textit{KB}}^\alpha$ is not the torus.  If $C_1$ is homologous to $C_1+C_2$, then it follows that \[\langle C_1, C_1 \rangle = \langle C_v+C_2, C_v+C_2 \rangle= \langle C_v+C_2,C_v+C_2\rangle = 1.\]  Thus, Remark \ref{remark:TransverseCrossings} implies that $GP(p,2)\mycolon C_1$ is an orientation-reversing circle, and by Lemma \ref{lemma:OrientationReversingVoltage}, $\mbox{\textit{KB}}^\alpha$ is not the torus.\medskip

\noindent\textit{Case 3d}: $GP(p,2)$ in the torus\medskip

\noindent\textit{Subcase 1}: $C_v$ is homologically trivial\medskip

If $C_v$ is homologically trivial, then, by Lemma \ref{lemma:TDecomposition}, $GP(p,2)\mycolon C_v$ is the boundary of the punctured torus and punctured sphere forming the connected sum $T=T \# S^2$.  Therefore, $GP(p,2)\mycolon C_v$ is contractible.  Since $GP(p,2)\setminus E(C_v)$ is connected, it follows that $GP(p,2)\setminus GP(p,2)\mycolon C_v$ is contained in either the punctured torus or the punctured sphere.  Since $GP(p,2)$ is nonplanar, it follows that $GP(p,2)\mycolon C_v$ bounds a face $f$, and so by Theorem \ref{theorem:RHEquation}, we have that 
\[\chi(T^\alpha)\le 2\cdot 0 -1 = -1,\]
which precludes $T^\alpha$ from being the torus.\medskip

\noindent\textit{Subcase 2}: $C_v$ is homologically nontrivial\medskip

\begin{figure}[H]
\begin{center}
\includegraphics[scale=.45]{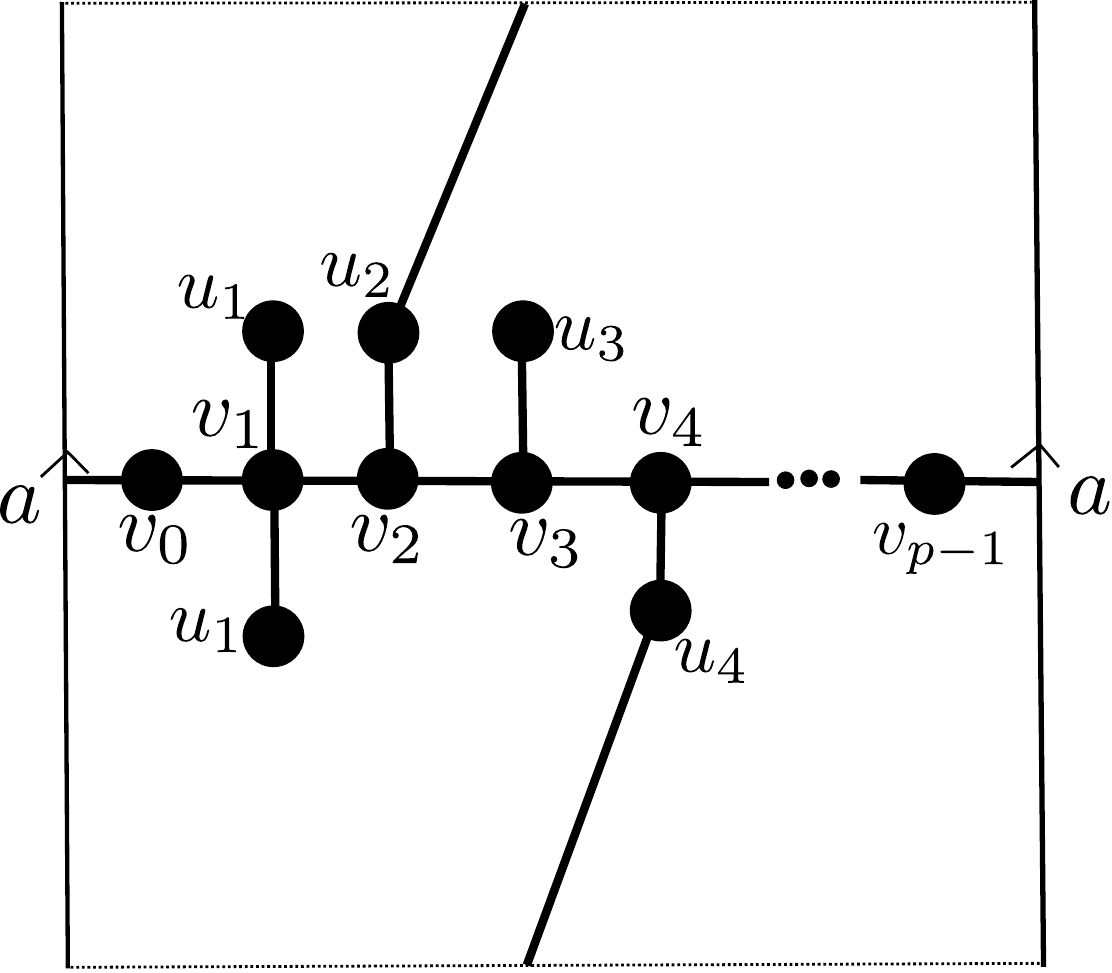}
\caption{A ribbon neighborhood of $GP(p,2)\mycolon C_v$ partially embedded in the torus, with the two possible placements for $(v_1,u_1)$.}\label{fig:GPp2TorusBase1}
\end{center}
\end{figure}

Since the torus is orientable, $GP(p,2)\mycolon C_v$ is an orientation-preserving circle and a ribbon neighborhood $R(C_v)$ is homeomorphic to an annulus.  Moreover Since $C_v$ is homologically nontrivial, $GP(p,2)\mycolon C_v$ is a nonseparating circle, and so, there is another element of our chosen bases for $Z(GP(p,2))$, say $C_2$, that satisfies $\langle C_v,C_2\rangle=1$.  Consider Figure \ref{fig:GPp2TorusBase1}, and note that without loss of generality, $(v_3,u_3)$ can be drawn as it is, above $GP(p,2)\mycolon C_v$.  Also note the two possible placements for the edge $(v_1,u_1)$ are also shown.

No matter how we draw the edge $(v_1,u_1)$, the fundamental cycles $C_1$ and $C_2$ must satisfy $\langle C_1,C_2\rangle=1$.  So, Lemma \ref{lemma:ZIntersections} implies that $C_1$ and $C_2$ are homologically nontrivial.  Since $GP(p,2)\mycolon C_1$ is orientation preserving, we have that
\[\langle C_1,C_1\rangle = \langle C_2,C_2\rangle = 0.\]Since there are two circles in $GP(p,2)^\alpha$ forming the fiber over $GP(p,2)\mycolon C_1$ and $GP(p,2)\mycolon C_2$, the fact that $p\colon T^\alpha \rightarrow T$ implies the following: 
\[\langle C_1^0,C_2^1\rangle = 1,\ \langle C_1^0, C_2^0\rangle = 0,\] and
\[\langle C_1^1, C_2^1\rangle = 1,\ \langle C_1^1,C_2^0\rangle = 0.\] We conclude that if we let $X=\left \{C_1^0,C_2^1,C_1^1,C_2^0\right \rbrace$, then  
\[M_{X_2}= \begin{pmatrix} 0 & 1 & 0 & 0 \\ 1 & 0 & 0 & 0 \\ 0 & 0 & 0 & 1 \\ 0 & 0 & 1 & 0 \end{pmatrix}.\]
Since the rows of $M_{X}$ are linearly independent, Theorem \ref{theorem:MatrixIndependence} implies that $\beta_1(T^\alpha)\ge 4$, which precludes the possibility that $\mbox{\textit{KB}}^\alpha$ is the torus since $\beta_1(T)=2$.\endproof

\noindent \textit{Proof of Part \ref{theorem:ThePointPart5}.}\medskip

\begin{figure}[H]
\begin{center}
\includegraphics[scale=.45]{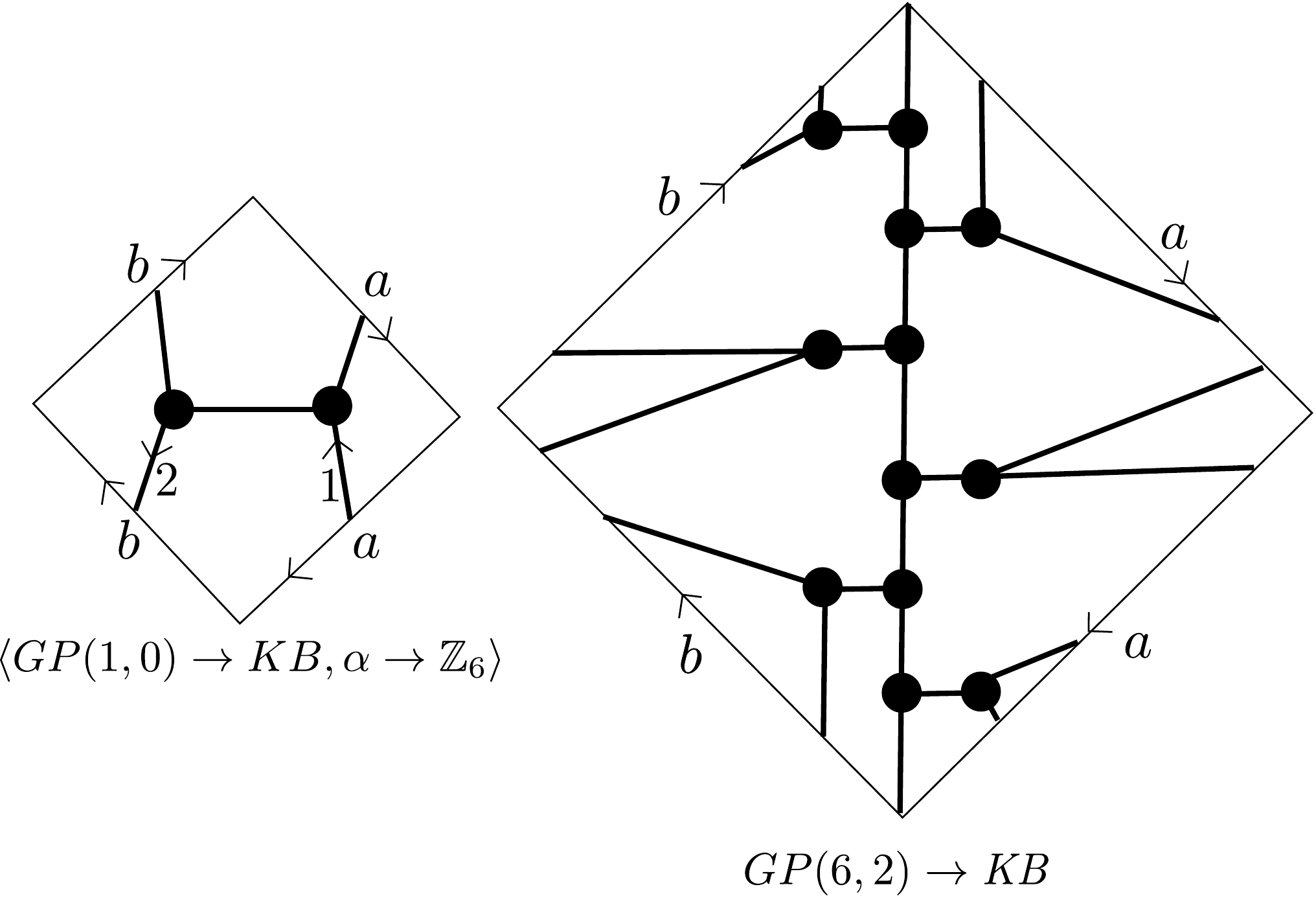}
\caption{An ordinary voltage graph embedding and its derived embedding, which is $GP(6,2)$ embedded in the Klein bottle.  The zero voltages are not shown in the base embedding.}\label{fig:KBBase1}
\end{center}
\end{figure}

Let $q$ be an odd prime.  Figure \ref{fig:KBBase1} verifies the case for which $q=3$.

\begin{figure}[H]
\begin{center}
\includegraphics[scale=.45]{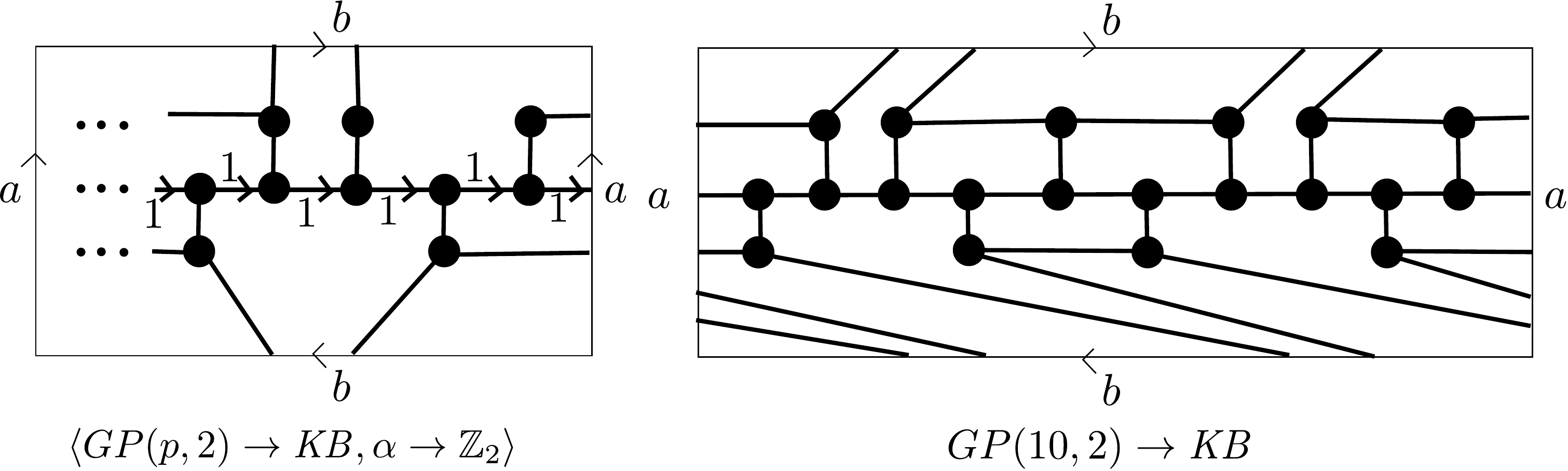}
\caption{An ordinary voltage graph embedding and a special case of its derived embedding, which is $GP(2p,2)\rightarrow KB$ for $p\ge5$.  The derived embedding is the special case for which $p=5$.  The zero voltages are not shown in the base embedding.}\label{fig:KBBase2}
\end{center}
\end{figure}

Figure \ref{fig:KBBase2} verifies the case for which $q\ge5$.

\endproof

\begin{remark}The case of Part \ref{theorem:BigCosetTheoremPart4} of Theorem \ref{theorem:ThePoint} for the case that $p=5$ is still open  To decide this case, one must know all of the free-actions of groups on $GP(10,2)$.  Since \cite[Theorem 1]{FGW} states that $GP(10,2)$ is vertex transitive, there are more actions to catalogue than those described in the proof of Part \ref{theorem:BigCosetTheoremPart4} of Theorem \ref{theorem:ThePoint}. \end{remark}

\section{Acknowledgements}
The author wishes to thank his dissertation advisor, Lowell Abrams, for his wisdom and counsel.  Most of the content of this article appears in the author's dissertation.

\bibliographystyle{plain}
\bibliography{OrdinaryVoltageGraphsAndDerivedHomology.bib}

\end{document}